\documentclass[twoside,11pt,a4paper,leqno]{article}
\usepackage{amssymb,amsmath,amscd,euscript,verbatim,array}
\usepackage[center,pagestyles]{titlesec}
\usepackage{geometry}
\usepackage{anysize}
\usepackage{fancyhdr}
\usepackage{indentfirst}
\usepackage{graphicx}
\usepackage{color}
\usepackage{ifpdf}
\usepackage{mathrsfs}
\usepackage{fancyref}
\usepackage{cite}
\usepackage[numbers,sort&compress]{natbib}
\marginsize{3.5cm}{3.5cm}{2cm}{2cm}

\titleformat{\section}{\centering\normalsize}{\thesection.}{0.5em}{}
\titleformat{\subsection}{\normalsize\bfseries}{\thesubsection.}{0.5em}{}
\titleformat{\subsubsection}{\normalsize\bfseries}{\thesubsubsection.}{0.5em}{}
\newcommand{\N}{\mathbb{N}}

\newcommand{\R}{\mathbb{R}}

\newtheorem{Theorem}{Theorem}[section]
\newtheorem{Definition}[Theorem]{Definition}
\newtheorem{Lemma}[Theorem]{Lemma}
\newtheorem{Exercise}[Theorem]{Exercise}
\newtheorem{Proposition}[Theorem]{Proposition}
\newtheorem{Remark}[Theorem]{Remark}

\newcommand{\gm}{\gamma}

\newcommand{\A}{\mathbb{A}}

\def\a{a}

\newcommand{\bthm}{\begin{Theorem}}
\newcommand{\ethm}{\end{Theorem}}
\newcommand{\bpr}{\begin{Proposition}}
\newcommand{\epr}{\end{Proposition}}
\newcommand{\blm}{\begin{Lemma}}
\newcommand{\elm}{\end{Lemma}}
\newcommand{\bex}{\begin{Exercise}}
\newcommand{\eex}{\end{Exercise}}
\newcommand{\be}{\begin{equation}}
\newcommand{\ee}{\end{equation}}
\newcommand{\beal}{\begin{aligned}}
\newcommand{\enal}{\end{aligned}}
\newcommand{\brm}{\begin{Remark}}
\newcommand{\erm}{\end{Remark}}
\newcounter{item}[section]

\newcommand{\Proof}{\textbf{Proof}\hspace{0.3cm}}
\newcommand{\End}{\ensuremath{\hfill{\Box}}\\}
\renewcommand{\title}[1]{\begin{center}\textbf{\large #1}\end{center}}
\renewcommand{\author}[1]{\begin{center}\normalsize #1\end{center}}
\renewcommand{\date}[1]{\begin{center}#1\end{center}}

\makeatletter \@addtoreset{equation}{section}
\makeatother
\begin{document}
\vspace{10pt}
\title{WEAK KAM THEORY FOR  GENERAL HAMILTON-JACOBI EQUATIONS I: THE SOLUTION SEMIGROUP UNDER PROPER CONDITIONS}

\vspace{6pt}
\author{\sc Xifeng Su, Lin Wang and Jun Yan}

\vspace{10pt} \thispagestyle{plain}

\begin{quote}
\small {\sc Abstract.}We consider the following evolutionary Hamilton-Jacobi equation with initial condition:
\begin{equation*}
\begin{cases}
\partial_tu(x,t)+H(x,u(x,t),\partial_xu(x,t))=0,\\
u(x,0)=\phi(x).
\end{cases}
\end{equation*}
Under some assumptions on  $H(x,u,p)$ with respect to $p$ and $u$, we provide a variational principle on the  evolutionary Hamilton-Jacobi equation.  By  introducing an implicitly defined  solution semigroup, we  extend Fathi's weak KAM theory to certain more general cases,  in which $H$ explicitly depends on the unknown function $u$. As an application, we show  the viscosity solution of the evolutionary Hamilton-Jacobi equation with initial condition tends asymptotically to the weak KAM solution of the following stationary Hamilton-Jacobi equation:
\begin{equation*}
H(x,u(x),\partial_xu(x))=0.
\end{equation*}.
\end{quote}
\begin{quote}
\end{quote}
\begin{quote}
\end{quote} \vspace{25pt}

\tableofcontents
\newpage
\section{\sc Introduction and main results}
Let $M$ be a compact connected  $C^2$ manifold and $H:T^*M\times\R\rightarrow\R$ be a $C^2$ function called a Hamiltonian. For a given $T>0$, we consider the following Hamilton-Jacobi equation:
\begin{equation}\label{hje}
\partial_tu(x,t)+H(x,u(x,t),\partial_xu(x,t))=0,
\end{equation}where $(x,t)\in M\times[0,T]$ and with the initial condition:
\[u(x,0)=\phi(x).\]
The characteristics of (\ref{hje}) satisfies the following equation:
\begin{equation}\label{hjech}
\begin{cases}
\dot{x}=\frac{\partial H}{\partial p},\\
\dot{p}=-\frac{\partial H}{\partial x}-\frac{\partial H}{\partial u}p,\\
\dot{u}=\frac{\partial H}{\partial p}p-H.
\end{cases}
\end{equation}

In 1983, M. Crandall and P. L. Lions introduced a notion of weak solutions of (\ref{hje}) named viscosity solution for overcoming the lack of uniqueness of the solution due to the crossing of characteristics (see \cite{Ar,CEL,CL2}).
During the same period, S. Aubry and J. Mather developed a seminal work so called Aubry-Mather theory on global action minimizing orbits for area-preserving twist maps (see \cite{Au,AD,M0,Mat,M2,M3} for instance). In 1991, J. Mather generalized the Aubry-Mather theory to  positive definite Lagrangian systems with multi-degrees of freedom (see \cite{M1}).

There is a close connection between viscosity solutions  and Aubry-Mather theory. Roughly speaking, the global minimizing orbits in Aubry-Mather theory can be embedded into the characteristic fields of PDEs. The similar ideas were reflected  in pioneering papers \cite{E} and \cite{F2} respectively. In \cite{E}, W. E  was concerned with certain weak solutions of the Burgers equation. In \cite{F2}, A. Fathi considered the Hamilton-Jacobi equations under so called Tonelli conditions (see (H1)-(H3) below). In particular, the Hamiltonian $H$ does not explicitly depend on $u$. He introduced a weak solution named weak KAM solution  and implied that the weak KAM solution is a viscosity solution, which initiated so called weak KAM theory.   A systematic introduction to weak KAM theory can be found in \cite{F3}.

In this paper, we are devoted to exploring the dynamics of more general Hamilton-Jacobi equations, in which $H$ explicitly depends on the unknown function $u$. Precisely speaking,
we are concerned with a   $C^2$ Hamiltonian $H(x,u,p)$ satisfying the following assumptions:
\begin{itemize}
\item [\textbf{(H1)}] \textbf{Positive Definiteness}: $H(x,u,p)$ is strictly convex with respect to $p$;
\item [\textbf{(H2)}] \textbf{Superlinear Growth}: For every compact set $I\subset \R$ and every $u\in I$, $H(x,u,p)$ is superlinear growth with respect  to $p$;
\item [\textbf{(H3)}] \textbf{Completeness of the Flow}: The flows of (\ref{hjech}) generated by $H(x,u,p)$ are complete.
\item [\textbf{(H4)}] \textbf{Uniform Lipschitz}: $H(x,u,p)$ is uniformly Lipschitz with respect to $u$.
\item [\textbf{(H5)}] \textbf{Monotonicity}: $H(x,u,p)$ is non-decreasing with respect to $u$.
\end{itemize}
Literately, (H1)-(H3) are called Tonelli conditions (see \cite{F3,M1}). (H5) is referred to as ``proper" condition (see \cite{CHL2}).

 Under the assumptions (H1)-(H5), we provide a variational principle on the  evolutionary Hamilton-Jacobi equation (\ref{hje}).  By  introducing an implicitly defined  solution semigroup, we  extend Fathi's weak KAM theory to certain more general cases, in which $H$ explicitly depends on $u$. As an application, we show  the viscosity solution of the evolutionary Hamilton-Jacobi equation with initial condition tends asymptotically to the weak KAM solution of the stationary Hamilton-Jacobi equation.

 The aim of this paper is to show the main ideas of exploring the dynamics of  more general Hamilton-Jacobi equations. To avoid the digression, we do not discuss whether  the assumptions (H1)-(H5) are optimal, which will be focused in the future work.

To state the main results, we first introduce some technology.
We use $\mathcal{L}: T^*M\rightarrow TM$ to denote  the Legendre transformation. Let
$\bar{\mathcal{L}}:=(\mathcal{L},  Id)$, where $Id$ denotes the identity map from $\R$ to $\R$. Then $\bar{\mathcal{L}}$ denotes a diffeomorphism from $T^*M\times\R$ to $TM\times\R$. By $\bar{\mathcal{L}}$,
the Lagrangian $L(x,u, \dot{x})$ associated to $H(x,u,p)$ can be denoted by
\[L(x,u, \dot{x}):=\sup_p\{\langle \dot{x},p\rangle-H(x,u,p)\}.\]
Let $\Psi_t$ denote the flow generated by $H(x,u,p)$. The flow generated by $L(x,u,\dot{x})$
can be denoted by $\Phi_t:=\bar{\mathcal{L}}\circ\Psi_t\circ\bar{\mathcal{L}}^{-1}$. From
(H1)-(H5), it follows that the Lagrangian $L(x,u,\dot{x})$ satisfies:
\begin{itemize}
\item [\textbf{(L1)}]  \textbf{Positive Definiteness}: $L(x,u,\dot{x})$ is strictly convex with respect  to $\dot{x}$;
\item [\textbf{(L2)}] \textbf{Superlinear Growth}: For every compact set $I\subset \R$ and every $u\in I$, $L(x,u,\dot{x})$  is uniformly superlinear growth with respect  to $\dot{x}$;
\item [\textbf{(L3)}] \textbf{Completeness of the Flow}: The flows generated by $L(x,u,\dot{x})$  are complete.
\item [\textbf{(L4)}] \textbf{Uniform Lipschitz}: $L(x,u,\dot{x})$ is uniformly Lipschitz with respect to $u$.
\item [\textbf{(L5)}] \textbf{Monotonicity}: $L(x,u,\dot{x})$ is non-increasing with respect to $u$.
\end{itemize}

 If a Hamiltonian $H(x,u,p)$ satisfies (H1)-(H5) (associated  $L(x,u,\dot{x})$ satisfying (L1)-(L5)), then we obtain the following theorem:
\begin{Theorem}\label{one}
 There exists a unique   $u(x,t)\in C(M\times[0,T],\R)$ satisfying $u(x,0)=\phi(x)$ such that
\begin{equation}\label{uxt}
u(x,t)=\inf_{\gm(t)=x}\left\{\phi(\gm(0))+\int_0^tL(\gm(\tau),u(\gm(\tau),\tau),\dot{\gm}(\tau))d\tau\right\}.
\end{equation}where the infimums are taken among the absolutely continuous curves $\gm:[0,t]\rightarrow M$. The infimums can be attained at a $C^1$ curve denoted by $\bar{\gm}$. Moreover, for $\tau\in (0,t)$,
$(\bar{\gm}(\tau),\bar{u}(\tau),p(\tau))$ satisfies the characteristics equation (\ref{hjech}) where
\begin{equation*}
\bar{u}(\tau)=u(\bar{\gm}(\tau),\tau)\quad\text{and}\quad p(\tau)=\frac{\partial L}{\partial \dot{x}}(\bar{\gm}(\tau),u(\bar{\gm}(\tau),\tau),\dot{\bar{\gm}}(\tau)).
\end{equation*}
\end{Theorem}

By analogy of the notion of weak KAM solution of the Hamilton-Jacobi equation independent of  $u$  (see \cite{F3}). We define another weak solution of (\ref{hje})  called variational solution (see Definition \ref{nw}). Based on Theorem \ref{one}, we construct a variational solution of (\ref{hje}) with initial condition.  Following \cite{F3}, we show that the variational solution of (\ref{hje}) is a viscosity solution of (\ref{hje}). Based on the uniqueness of the viscosity solution under (H1)-(H5), we have the following theorem.
\begin{Theorem}\label{two}
For any $\phi(x)\in C(M,\R)$, there exists a unique viscosity solution $u(x,t)$ of (\ref{hje}) with initial condition $u(x,0)=\phi(x)$.
\end{Theorem}

Theorem \ref{one} provides a variational principle on the  evolutionary Hamilton-Jacobi equation as (\ref{hje}), from which
there exists an implicitly defined semigroup denoted by $T_t$ such that $u(x,t)=T_t\phi(x)$, where $u(x,t)$ satisfies (\ref{uxt}). To fix the notion, we call $T_t$ a solution semigroup.
We use $c(L(x,a,\dot{x}))$ to denote the Ma\~{n}\'{e} critical value of $L(x,a,\dot{x})$. By \cite{CIPP}, we have
\begin{equation}
c(L(x,a,\dot{x}))=\inf_{u\in C^1(M,\R)}\sup_{x\in M}H(x,a,\partial_xu).
\end{equation}
Let
\begin{equation}
\mathcal{C}=\left\{c(L(x,a,\dot{x})): a\in\R\right\}.
\end{equation}It follows that for any $c\in \mathcal{C}$, there exists $a\in\R$ such that $c(L(x,a,\dot{x}))=c$. Let
$L_c=L+c$, then $c(L_c(x,a,\dot{x}))=0$. Without ambiguity, we still use $L$ instead of $L_c$. The same to $H$ and $T_t$.
Let $\|\cdot\|_\infty$ be $C^0$-norm. We have the following theorem:
\begin{Theorem}\label{four}
For any $\phi(x),\psi(x)\in C(M,\R)$ and $t\geq 0$, the solution semigroup $T_t$ has following properties:
\begin{itemize}
\item [I.]  for $\phi\leq\psi$, $T_t\phi\leq T_t\psi$;
\item [II.] $\|T_t\phi-T_t\psi\|_\infty\leq \|\phi-\psi\|_\infty$;
\item [III.] there exists a positive constant $K$ independent of $t$ such that $\|T_t\phi\|_\infty\leq K$;
\item [IV.] for  $\delta>0$, the family of functions $(x,t)\rightarrow T_t\phi(x)$ is equi-Lipschitz on $(x,t)\in M\times [\delta,+\infty)$.
\end{itemize}
\end{Theorem}

For the autonomous  systems with Lagrangian $L(x,\dot{x})$, the convergence of so called Lax-Oleinik semigroup was established in \cite{F22}. By \cite{FaM}, such convergence fails for the non-autonomous Lagrangian systems. A new kind of operators was found in \cite{WY} for the time periodic Lagrangian systems to overcome the failure of the convergence of the Lax-Oleinik semigroup. Different from the previous results, the solution semigroup $T_t$ we consider is associated to $L(x,u(x,t),\dot{x})$, which is neither autonomous nor periodic with respect to $t$. It results in the lack of conservation of energy of the system and compactness of the underlying manifold. Hence, it is necessary  for establishing the convergence of $T_t$ to find a new way.
Based on Theorem \ref{one} and Theorem \ref{four}, we  obtain the convergence of the solution semigroup $T_t$ by considering the evolution of $H$ along the characteristics.
\begin{Theorem}\label{five}
For any $\phi(x)\in C(M,\R)$, $\lim_{t\rightarrow\infty}T_t\phi(x)$ exists. Moreover, let \[u_\infty(x):=\lim_{t\rightarrow\infty}T_t\phi(x),\]
 then $u_\infty$ is a weak KAM solution of the following stationary equation:
\begin{equation*}
H(x, u(x), \partial_x u(x))=0.
\end{equation*}
\end{Theorem}

By inspiration of \cite{F22}, the large time behavior of viscosity solutions of Hamilton-Jacobi equations with Hamiltonian independent of $u$ was explored comprehensively  based on both dynamical  and PDE approaches
 (see \cite{ds, II, NR} for instance). Theorem \ref{two} implies $u(x,t):=T_t\phi(x)$ is the unique viscosity solution of (\ref{hje}) with initial condition $u(x,0)=\phi(x)$.
As an application of Theorem \ref{five}, we  obtain the large time behavior of viscosity solutions of (\ref{hje}). More precisely, we have
\begin{Theorem}\label{six}
Let $u(x,t)$ be a viscosity solution of (\ref{hje}), then $\lim_{t\rightarrow \infty}u(x,t)$ is a weak KAM solution of the following stationary equation:
\begin{equation*}
H(x, u(x), \partial_x u(x))=0.
\end{equation*}
\end{Theorem}

 This paper is outlined as follows. In Section 2, some definitions are recalled as preliminaries. In Section 3, an implicitly variational principle is established.  Moreover, Theorem \ref{one} can be obtained. In Section 4, a construction of the viscosity solution is provided, which implies Theorem \ref{two}.  In Section 5, an implicitly defined solution semigroup is introduced and some properties are detected, from which Theorem \ref{four} is proved.  In Section 6, following from Theorem \ref{four}, the convergence of the solution semigroup is shown. Moreover, both of Theorem \ref{five} and Theorem \ref{six} can be verified.

\section{\sc Preliminaries}
In this section, we recall the definitions of the weak KAM solution and the viscosity solution of (\ref{hje})  (see \cite{CEL,CL2,F3}) and some aspects of Mather-Fathi theory for  the sake of completeness.
\subsection{Weak KAM solutions and viscosity solutions}
A function $H(x,p):TM\rightarrow\R$ called a Tonelli Hamiltonian if $H$ satisfies (H1)-(H2). For the autonomous Hamiltonian systems, the assumption (H3) holds obviously  from the compactness of $M$. The associated Lagrangian is denoted by $L(x,\dot{x})$ via the Legendre transformation. In \cite{F1}, Fathi introduced the definition of the weak KAM solution of negative type of the following Hamilton-Jacobi equation:
\begin{equation}\label{fathi}
H(x,  \partial_x u(x))=0,\quad x\in M,
\end{equation}where $H$ is a Tonelli Hamiltonian.
\begin{Definition}\label{weakkk}
A function $u\in C(M,\R)$ is called a  weak KAM solution of negative type of (\ref{fathi}) if
\begin{itemize}
\item [(i)] for each continuous piecewise $C^1$ curve $\gm:[t_1,t_2]\rightarrow M$ where $t_2>t_1$, we have
\begin{equation}
u(\gm(t_2))-u(\gm(t_1))\leq\int_{t_1}^{t_2}L(\gm(\tau),\dot{\gm}(\tau))d\tau;
\end{equation}
\item [(ii)] for any $x\in M$, there exists a $C^1$ curve $\gm:(-\infty,0]\rightarrow M$ with $\gm(0)=x$ such that for any $t\in (-\infty,0]$, we have
  \begin{equation}
u(x)-u(\gm(t))=\int_{t}^{0}L(\gm(\tau),\dot{\gm}(\tau))d\tau.
\end{equation}
\end{itemize}
\end{Definition}
By  analogy of the definition above, it is easy to define the weak KAM solution of negative type of more general Hamilton-Jacobi equation as follows:
\begin{equation}\label{fathiw}
H(x, u(x), \partial_x u(x))=0,\quad x\in M.
\end{equation}
\begin{Definition}\label{weakkam}
A function $u\in C(M,\R)$ is called a  weak KAM solution of negative type of (\ref{fathiw}) if
\begin{itemize}
\item [(i)] for each continuous piecewise $C^1$ curve $\gm:[t_1,t_2]\rightarrow M$ where $t_2>t_1$, we have
\begin{equation}
u(\gm(t_2))-u(\gm(t_1))\leq\int_{t_1}^{t_2}L(\gm(\tau),u(\gm(\tau)),\dot{\gm}(\tau))d\tau;
\end{equation}
\item [(ii)] for any $x\in M$, there exists a $C^1$ curve $\gm:(-\infty,0]\rightarrow M$ with $\gm(0)=x$ such that for any $t\in (-\infty,0]$, we have
  \begin{equation}
u(x)-u(\gm(t))=\int_{t}^{0}L(\gm(\tau),u(\gm(\tau)),\dot{\gm}(\tau))d\tau.
\end{equation}
\end{itemize}
\end{Definition}

Following from \cite{CEL,CL2,F3},
 a viscosity solution of (\ref{hje}) can be defined as follows:
\begin{Definition}\label{visco}
Let $V$ be an open subset  $V\subset M$,
\begin{itemize}
\item [(i)] A function $u:V\times[0,T]\rightarrow \R$ is a subsolution of  (\ref{hje}), if for every $C^1$ function $\phi:V\times[0,T]\rightarrow\R$ and every point $(x_0,t_0)\in V\times[0,T]$ such that $u-\phi$ has a maximum at $(x_0,t_0)$, we have
\begin{equation}
\partial_t\phi(x_0,t_0)+H(x_0,t_0,u(x_0,t_0),\partial_x\phi(x_0,t_0))\leq 0;
\end{equation}
\item [(ii)] A function $u:V\times[0,T]\rightarrow \R$ is a supersolution of  (\ref{hje}), if for every $C^1$ function $\psi:V\times[0,T]\rightarrow\R$ and every point $(x_0,t_0)\in V\times[0,T]$ such that $u-\psi$ has a minimum at $(x_0,t_0)$, we have
\begin{equation}
\partial_t\psi(x_0,t_0)+H(x_0,t_0,u(x_0,t_0),\partial_x\psi(x_0,t_0))\geq 0;
\end{equation}
\item [(iii)] A function $u:V\times[0,T]\rightarrow \R$ is a viscosity solution of  (\ref{hje}) if it is both a subsolution and a supersolution.
\end{itemize}
\end{Definition}
Under the assumptions (H1)-(H5), it follows from the comparison theorem that the viscosity solution of (\ref{hje}) with initial condition is unique (see \cite{CL2}) and it is a locally semiconcave function (see \cite{CS}).

Both of Definition \ref{weakkk} and Definition \ref{weakkam} are concerned with the weak KAM solutions defined on $M\times\R$, while the viscosity solutions of (\ref{hje}) are defined on $M\times[0,T]$. As a bridge connecting them,
we give the definition of another weak solution of (\ref{hje}) with initial condition called a variational solution.
\begin{Definition}\label{nw}
For a given $T>0$, a variational solution of (\ref{hje}) with initial condition is a function $u: M\times [0,T]\rightarrow\R$ for which the following are satisfied:
\begin{itemize}
\item [(i)] for each continuous piecewise $C^1$ curve $\gm:[t_1,t_2]\rightarrow M$ where $0\leq t_1<t_2\leq T$, we have
\begin{equation}
u(\gm(t_2),t_2)-u(\gm(t_1),t_1)\leq\int_{t_1}^{t_2}L(\gm(\tau),u(\gm(\tau),\tau),\dot{\gm}(\tau))d\tau;
\end{equation}
\item [(ii)] for any $0\leq t_1<t_2\leq T$ and $x\in M$, there exists a $C^1$ curve $\gm:[t_1,t_2]\rightarrow M$ with $\gm(t_2)=x$ such that
  \begin{equation}
u(x,t_2)-u(\gm(t_1),t_1)=\int_{t_1}^{t_2}L(\gm(\tau),u(\gm(\tau),\tau),\dot{\gm}(\tau))d\tau.
\end{equation}
\end{itemize}
\end{Definition}
The existence of  the variational solutions will be verified in Section 4.

\subsection{The minimal action and the Peierls barrier}
Let $L:TM\rightarrow\R$ be a Tonelli Lagrangian. We define the function $h_t:M\times M\rightarrow\R$ by
\begin{equation}\label{mather}
h_t(x,y)=\inf_{\substack{\gm(0)=x\\ \gm(t)=y}}\int_0^t L(\gm(\tau),\dot{\gm}(\tau))d\tau,
\end{equation}
where the infimums are taken among the absolutely continuous curves $\gm:[0,t]\rightarrow M$.
By Tonelli theorem (see \cite{F3,M1}), the infimums in (\ref{mather}) can be achived. Let $\bar{\gm}$ be an absolutely continuous curve with $\bar{\gm}(0)=x$ and $\bar{\gm}(t)=y$ such that the infinmum is achieved at $\bar{\gm}$. Then $\bar{\gm}$ is called a minimal curve. By \cite{M1}, the minimal curves satisfy the Euler-Lagrange equation generated by $L$. The quantity $h_t(x,y)$ is called a minimal  action.

From the definition of $h_t(x,y)$, it follows that for each $x,y,z\in M$ and each $t,t'>0$, we have
 \begin{equation}\label{yan}
h_{t+t'}(x,z)\leq h_{t}(x,y)+h_{t'}(y,z).
\end{equation}
  In particular,  we have
\begin{equation}\label{wang}
h_{t+t'}(x,y)=h_{t}(x,\bar{\gm}(t))+h_{t'}(\bar{\gm}(t),y),
\end{equation}
where $\bar{\gm}$ is  a minimal curve with $\bar{\gm}(0)=x$ and $\bar{\gm}(t+t')=y$.

We use $c[0]$ to denote the Ma\~{n}\'{e} critical value of $L(x,\dot{x})$. By \cite{CIPP}, we have
\begin{equation}
c[0]=\inf_{u\in C^1(M,\R)}\sup_{x\in M}H(x,\partial_xu).
\end{equation}
The definition of the Peierls barrier is due to Mather (see \cite{Mvc}). A systematic introduction can be founded in \cite{F3}.
The Peierls barrier  is the function $h:M\times M\rightarrow\R$ defined by
\[h(x,y)=\liminf_{t\rightarrow\infty}h_t(x,y)+c[0]t.\]
For the autonomous Lagrangian, ``liminf" can be replaced by ``lim". By the compactness of $M$, for each $t_0>0$, there exists a constant $C_{t_0}$ such that for any $t\geq t_0$ and $x,y\in M$,
\begin{equation}\label{lin}
|h_t(x,y)+c[0]t|\leq C_{t_0},
\end{equation}which implies the values of the map $h$ is finite (see Lemma 5.3.2 in \cite{F3}). The property above play a crucial role in our context.

\section{\sc Variational principle}
For every given continuous function $\phi$ on $M$, we now define
the operator $\A:C(M\times [0,T], \mathbb{R}
)\circlearrowleft$ depending on $\phi$ as follows:
\begin{equation}\label{operator A}
\begin{split}
\A[u](x,t) = \inf_{\substack{\gamma(t)=x\\ \gamma\in
C^{ac}([0,t],M)}} \left\{ \phi\big(\gamma(0)\big) + \int_0^t
L\big(\gamma(s),u\big(\gamma(s),s\big), \dot\gamma(s)\big) \ ds
\right\},
\end{split}
\end{equation}where $u\in C(M\times[0,T], \mathbb{R})$.
By means of a simple modification of Tonelli's theorem (see \cite{F3} and \cite{M1}), we have the following lemma.
\begin{Lemma}\label{tonelli}
For a given $u\in C(M\times[0,T],\R)$, the infimums in (\ref{operator A}) are attained at  absolutely continuous curves with the end point $x$.
\end{Lemma}
We omit the proof of Lemma \ref{tonelli} for the consistency of the context.
\subsection{The fixed point of $\A$}

In the following, we will prove that the operator $\A$ has a unique fixed point.
\begin{Lemma}\label{fixed point existence}
$\A$ has a unique fixed point.
\end{Lemma}
\Proof
From (L5), it follows that for  every $u,v\in C(M\times[0,T], \mathbb{R})$,
\begin{equation*}
|L(x,u,\dot{x})-L(x,v,\dot{x})| \leq \lambda \ |u-v|,
\end{equation*}where $\lambda$ is a positive constant independent of $x$ and $\dot{x}$. Hence, for any given $t\in [0,T]$,
it follows from Lemma \ref{tonelli} that
\begin{equation*}
\begin{split}
&\big(\A[u] - \A[v]\big)(x,t)\\
\leq& \int_0^t (L\big(\gamma(s),u\big(\gamma(s),s\big),
\dot\gamma(s)\big) - L\big(\gamma(s),v\big(\gamma(s),s\big),
\dot\gamma(s)\big) )\
ds \\
\leq& \lambda \ \|u-v\|_{\infty} t
\end{split}
\end{equation*}where $\gamma\in C^{ac}([0,t],M)$ such that
\begin{equation}\label{u}
\A[v](x,t) = \phi\big(\gamma(0)\big) + \int_0^t
L\big(\gamma(s),v\big(\gamma(s),s\big), \dot\gamma(s)\big) \ ds.
\end{equation}

By exchange the position of $u$ and $v$, we obtain
\begin{equation*}
|\big(\A[u] - \A[v]\big)(x,t)| \leq \lambda \ \|u-v\|_{\infty}t.
\end{equation*}

Therefore, we have the following estimates:
\begin{equation*}
\begin{split}
 &\left|\big(\A^2[u] - \A^2[v]\big)(x,t)\right| \\
 \leq & \left| \int_0^t \lambda \big[ \A[u]\big(\gamma(s),s\big) -
 \A[v]\big(\gamma(s),s\big)\big] \ ds\right|\\
 \leq & \int_0^t  s \lambda^2 \|u-v\|_{\infty} \ ds \leq \frac{(t\lambda)^2}{2}
 \|u-v\|_{\infty}.
\end{split}
\end{equation*}

Moreover, continuing the above procedure, we obtain
\begin{equation}
\left|\big(\A^n[u] - \A^n[v]\big)(x,t)\right| \leq
\frac{(t\lambda)^n}{n!} \|u-v\|_{\infty}.
\end{equation}
Therefore, for any $t\in [0,T]$, there exists  $N\in \N$ large
enough such that $\A^{N}$ is a
contraction mapping and has a fixed point. That is, for any
$t\in [0,T]$ and $N\in\mathbb{N}$ large enough, there exists a
$u\in C(M, \mathbb{R})$ such that
\begin{equation}
\A^N[u](x) = u(x).
\end{equation}

We now show that $u$ is a fixed point of $\A$. Since
\begin{equation*}
\A[u] = \A\circ \A^N[u] = \A^N\circ \A[u],
\end{equation*}$\A[u]$ is also a fixed point of $\A^N$.
By the uniqueness of fixed point of contraction mapping, we have
\[A[u] = u.\]
This completes the proof of Lemma \ref{fixed point existence}.
\End

Lemma \ref{fixed point existence} shows that
 there exists $u(x,t)\in C(M\times [0,T],\R)$ such that
\begin{equation}\label{fixu}
u(x,t)=\inf_{\gm(t)=x}\left\{\phi(\gm(0))+\int_0^tL(\gm(\tau),u(\gm(\tau),\tau),\dot{\gm}(\tau))d\tau\right\}.
\end{equation}To fix the notions, we give a definition as follows:
\begin{Definition}\label{cali}
For $u(x,t)\in C(M\times[0,T],\R)$ satisfying (\ref{fixu}), a curve $\gm:I\rightarrow M$ is called a calibrated curve of $u$ if for andy $0\leq t_1< t_2\leq T$, we have
\[u(\gm(t_2),t_2)=u(\gm(t_1),t_1)+\int_{t_1}^{t_2}L(\gm(\tau),\tau,u(\gm(\tau),\tau),\dot{\gm}(\tau))d\tau.\]
\end{Definition}
\subsection{Calibrated curves and characteristics}
 In the following, we will show the relation between calibrated curves and characteristics of (\ref{hje}). More precisely, we have the following lemma:

\begin{Lemma}\label{chc}
Let $\bar{\gm}:[0,t]\rightarrow M$ be a calibrated curve of $u$, then $\bar{\gm}$ is $C^1$ and for $\tau\in (0,t)$,
$(\bar{\gm}(\tau),\bar{u}(\tau),p(\tau))$ satisfies the characteristics equation (\ref{hjech}) where
\begin{equation}\label{ccure}
\bar{u}(\tau)=u(\bar{\gm}(\tau),\tau)\quad\text{and}\quad p(\tau)=\frac{\partial L}{\partial \dot{x}}(\bar{\gm}(\tau),u(\bar{\gm}(\tau),\tau),\dot{\bar{\gm}}(\tau)).
\end{equation}
\end{Lemma}
\Proof
Since $\bar{\gm}\in C^{ac}([0,t],M)$, then the derivative $\dot{\bar{\gm}}(\tau)$ exists almost everywhere for $\tau\in [0,t]$. Let $t_0\in (0,t)$ be a differentiate point of $\bar{\gm}(\tau)$. For the simplicity of notations and without ambiguity, we denote
\begin{equation}
(x_0,u_0,v_0):=(\bar{\gm}(t_0),u(\bar{\gm}(t_0),t_0),\dot{\bar{\gm}}(t_0)).
 \end{equation}

 First of all, we will construct a classical solution on a cone-like region (see (\ref{tri}) below). Let $k:=|v_0|$ and
 \[B(0,2k):=\{v:|v|< 2k,\ v\in T_{x_0}M\}.\]
We use $B^*(0,2k)$ to denote the image of $B(0,2k)$ via the Legendre transformation $\mathcal{L}^{-1}:TM\rightarrow T^*M$. That is
\[B^*(0,2k):=\left\{p: p=\frac{\partial L}{\partial v}(x_0,u_0,v),\ v\in B(0,2k)\right\}.\]
  Let $\Psi_t:T^*M\times\R\rightarrow T^*M\times\R$ denote the follow generated by the characteristics equation (\ref{hjech}). Let $\pi$ be a projection from $T^*M\times\R$ to $ T^*M$ via $(x,p,u)\rightarrow (x,p)$  and let
$B_t^*(0,2k):=\pi\circ\Psi_{t-t_0}(B^*(0,2k),u_0)$. We denote
\[\Pi_t:B_t^*(0,2k)\rightarrow M.\]
Since the Legendre transformation $\mathcal{L}$ is a diffeomorphism, then for a given   $\epsilon>0$ small enough and $\tau\in [t_0,t_0+\epsilon]$,  $\Pi_\tau$ is a diffeomorphism onto the image denoted by $\Omega_\tau:=\Pi_\tau( B_\tau^*(0,2k))$.
 We use $\Omega^\epsilon$ to denote the following cone-like region:
\begin{equation}\label{tri}
\Omega^\epsilon:=\{(\tau,x): \tau\in (t_0,t_0+\epsilon),\ x\in \Omega_\tau\}.
 \end{equation}Then for any $(\tau,x)\in\Omega^\epsilon$, there exists a unique $p_0\in B^*(0,2k)$ such that $X(\tau)=x$ where
 \[(X(\tau),U(\tau),P(\tau)):=\Phi_\tau(x_0,u_0,p_0).\]
 Hence, for  any $(\tau,x)\in \Omega^\epsilon$, one can define a $C^1$ function by $S(x,\tau)=U(\tau)$. In particular, we have $S(x,t_0)=u_0$. Moreover, it follows from  the method of characteristics (see \cite{F3,Li} for instance) that $S(x,\tau)$
is a solution of the following equation:
\begin{equation}\label{sss}
\partial_\tau S(x,\tau)+H(x,S(x,\tau),\partial_x S(x,\tau))=0
\end{equation}with $\partial_x S(x_0,t_0)=\partial_v L(x_0,S(x_0,t_0),v_0)$, where $L$ denotes Lagrangian via the Legendre transformation associated to the Hamiltonian $H$. Fix $\tau\in [t_0,t_0+\epsilon]$ and let $S_\tau(x):=S(x,\tau)$. We denote
\begin{equation}
\text{grad}_LS_\tau(x):=\frac{\partial H}{\partial p}(x,S_\tau(x),p),
\end{equation}where $p=\partial_x S_\tau(x)$.  It is easy to see that $\text{grad}_LS_\tau(x)$ gives rise to a $C^1$ vector field on $M$. Moreover, we have the following claim.

\textbf{Claim:}  Let $\gm$ be an absolutely continuous curve with $(\tau,\gm(\tau))\in \Omega^\epsilon$ for  $\tau\in [a,b]\subset [t_0,t_0+\epsilon]$, we have
\begin{equation}\label{hjee}
S(\gm(b),b)-S(\gm(a),a)\leq\int_a^bL(\gm(\tau),
S(\gm(\tau),\tau),\dot{\gm}(\tau))d\tau,
\end{equation}
where the equality holds if and only if $\gm$ is a trajectory of the vector field $\text{grad}_LS_\tau(x)$.

\Proof
Since $S(x,\tau)$ is a $C^1$ function, then we have
\begin{equation}\label{sasb}
S(\gm(b),b)-S(\gm(a),a)=\int_a^b\left\{\frac{\partial S}{\partial t}(\gm(\tau),\tau)+\langle \frac{\partial S}{\partial x}(\gm(\tau),\tau),\dot{\gm}(\tau)\rangle\right\}d\tau.
\end{equation}By virtue of Fenchel inequality, for each $\tau$ where $\dot{\gm}(\tau)$ exists, we have
\begin{align*}
\langle \frac{\partial S}{\partial x}(\gm(\tau),\tau),\dot{\gm}(\tau)\rangle\leq &H(\gm(\tau),S(\gm(\tau),\tau),\frac{\partial S}{\partial x}(\gm(\tau),\tau))\\
&+L(\gm(\tau),
S(\gm(\tau),\tau),\dot{\gm}(\tau)).
\end{align*}It follows from (\ref{sss}) that for almost every $\tau\in [a,b]$
\begin{equation}
\frac{\partial S}{\partial t}(\gm(\tau),\tau)+\langle \frac{\partial S}{\partial x}(\gm(\tau),\tau),\dot{\gm}(\tau)\rangle\leq L(\gm(\tau),
S(\gm(\tau),\tau),\dot{\gm}(\tau)).
\end{equation}
By integration, it follows from (\ref{sasb})
\begin{equation}\label{ss}
S(\gm(b),b)-S(\gm(a),a)\leq\int_a^b L(\gm(\tau),
S(\gm(\tau),\tau),\dot{\gm}(\tau))d\tau.
\end{equation}We have equality in (\ref{ss}) if and only if the equality holds in the Fenchel inequality, i.e. $\dot{\gm}(\tau)=\text{grad}_LS_\tau(x)$ which means that $\gm$ is a trajectory of the vector field $\text{grad}_LS_\tau(x)$.
\End

The claim implies
\begin{equation*}
S(x,\tau)= u_0 + \inf_{\substack{\gm(t)=x \\  \gm(t_0)=x_0 } }\int_{t_0}^\tau L(\gm(s), S(\gm(s),s),\dot{\gm}(s))ds.
\end{equation*}
We recall $\bar{\gm}$ is a calibrated curve of $u$ with $\bar{\gm}(t_0)=x_0$. Since $\bar{\gm}$ is differentiable at $t_0$ with $|\dot{\bar{\gm}}(t_0)|=k$, then it follows from  the construction of $\Omega^\epsilon$ that for $\epsilon$ small enough and  $\tau\in [t_0,t_0+\epsilon]$,  we have $(\tau,\bar{\gm}(\tau))\in \Omega^\epsilon$.

\textbf{Claim:} For any $\tau\in [t_0,t_0+\epsilon]$, we have
\begin{equation}
S(\bar{\gm}(\tau),\tau)=u(\bar{\gm}(\tau),\tau).
\end{equation}

\Proof
By contradiction, we assume there exists $\tilde{t}\in [t_0,t_0+\epsilon]$ such that
\begin{equation}
S(\bar{\gm}(\tilde{t}),\tilde{t})\neq u(\bar{\gm}(\tilde{t}),\tilde{t}).
\end{equation}We only consider the case with $S(\bar{\gm}(\tilde{t}),\tilde{t})< u(\bar{\gm}(\tilde{t}),\tilde{t})$, the other case is similar. Let $\tilde{x}:=\bar{\gm}(\tilde{t})$ and let $\tilde{\gm}$ be a calibrated curve of $S$ with $\tilde{\gm}(t_0)=x_0$ and $\tilde{\gm}(\tilde{t})=\tilde{x}$.
We denote
\begin{equation}
F(\tau)=S(\tilde{\gm}(\tau),\tau)- u(\tilde{\gm}(\tau),\tau).
\end{equation}Hence, $F(\tau)$ is continuous and $F(t_0)=0$, $F(\tilde{t})<0$. Moreover, there exists $t_1\in [t_0,\tilde{t})$ such that $F(t_1)=0$ and $F(\tau)< 0$ for any $\tau\in (t_1,\tilde{t}]$, i.e.
\begin{equation}
S(\tilde{\gm}(\tau),\tau)< u(\tilde{\gm}(\tau),\tau).
\end{equation}By (L5), a simple calculation implies
\begin{align*}
&S(\tilde{x},\tilde{t})-u(\tilde{x},\tilde{t})\\
&\geq \int_{t_1}^{\tilde{t}}L(\tilde{\gm}(\tau),S(\tilde{\gm}(\tau),\tau),\dot{\tilde{\gm}}(\tau))
-L(\tilde{\gm}(\tau),u(\tilde{\gm}(\tau),\tau),\dot{\tilde{\gm}}(\tau))d\tau,\\
&\geq 0,
\end{align*}which is in contradiction with the assumption $S(\tilde{x},\tilde{t})<u(\tilde{x},\tilde{t})$. Hence,  for any $\tau\in [t_0,t_0+\epsilon]$, we have
\begin{equation}
S(\bar{\gm}(\tau),\tau)\geq u(\bar{\gm}(\tau),\tau).
\end{equation}Similarly, we have
\begin{equation}
S(\bar{\gm}(\tau),\tau)\leq u(\bar{\gm}(\tau),\tau).
\end{equation}Therefore,
\begin{equation}
S(\bar{\gm}(\tau),\tau)=u(\bar{\gm}(\tau),\tau),
\end{equation}which verifies the claim.
\End

 From the definition of $u$ (see (\ref{fixu})), it follows that
\begin{align*}
S(\bar{\gm}(t_0+\epsilon),t_0+\epsilon)=S(\bar{\gm}(t_0),t_0)
+\int_{t_0}^{t_0+\epsilon}L(\bar{\gm}(\tau),S(\bar{\gm}(\tau),\tau),\dot{\bar{\gm}}(\tau))d\tau,
\end{align*}which implies $\bar{\gm}(\tau)$  is a solution of the vector field $\text{grad}_LS_\tau(x)$. Let
\begin{equation}
\bar{u}(\tau):=u(\bar{\gm}(\tau),\tau)\quad\text{and}\quad p(\tau):=\frac{\partial L}{\partial \dot{x}}(\bar{\gm}(\tau),u(\bar{\gm}(\tau),\tau),\dot{\bar{\gm}}(\tau)).
\end{equation}Then $(\bar{\gm}(\tau),\bar{u}(\tau),p(\tau))$  is $C^1$ and satisfies the characteristics equation (\ref{hjech}).

 By (L3), a standard argument (see \cite{F3,M1}) shows that the differentiability of $\bar{\gm}(\tau)$ for $\tau\in [t_0,t_0+\epsilon]$ can be extended to the whole interval $(0,t)$. So far, we complete the proof of Lemma \ref{chc}.\End

So far, we complete the proof of Theorem \ref{one}.

\section{\sc Construction of  the viscosity solution}
In this section, we will provide a construction of the viscosity solution of (\ref{hje}). By Theorem \ref{one},  there exists a unique   $u(x,t)\in C(M\times[0,+\infty),\R)$ satisfying $u(x,0)=\phi(x)$ such that
\begin{equation}\label{fixu188}
u(x,t)=\inf_{\gm(t)=x}\left\{\phi(\gm(0))+\int_0^tL(\gm(\tau),u(\gm(\tau),\tau),\dot{\gm}(\tau))d\tau\right\}.
\end{equation}where the infimums are taken among the absolutely continuous curves. In particular, the infimums are attained at the characteristics of (\ref{hje}).

\begin{Lemma}\label{uisw}
$u(x,t)$ determined by (\ref{fixu188}) is a variational solution of (\ref{hje}) with initial condition $u(x,0)=\phi(x)$.
\end{Lemma}

\Proof
Let $\gm:[t_1,t_2]\rightarrow M$ be a continuous piecewise $C^1$ curve and Let $\bar{\gm}:[0,t_1]\rightarrow M$ be a calibrated curve of $u$ satisfying $\bar{\gm}(t_1)=\gm(t_1)$. We construct a curve $\xi:[0,t_2]\rightarrow M$ defined as follows:
\begin{equation}\label{xi}
\xi(t)=\left\{\begin{array}{ll}
\hspace{-0.4em}\bar{\gm}(t),& t\in [0,t_1],\\
\hspace{-0.4em}\gm(t),&t\in (t_1,t_2].\\
\end{array}\right.
\end{equation}
From (\ref{fixu}), it follows that
\begin{align*}
u(\gm(t_2)&,t_2)-u(\gm(t_1),t_1)\\
&= \inf_{\gm_2(t_2)=\gm(t_2)}\left\{\phi(\gm_2(0))+
\int_0^{t_2}L(\gm_2(\tau),u(\gm_2(\tau),\tau),\dot{\gm}_2(\tau))d\tau\right\}\\
&\ \ \ \ -\inf_{\gm_1(t_1)=\gm(t_1)}\left\{\phi(\gm_1(0))+
\int_0^{t_1}L(\gm_1(\tau),u(\gm_1(\tau),\tau),\dot{\gm}_1(\tau))d\tau\right\},\\
&\leq \phi(\xi(0))+
\int_0^{t_2}L(\xi(\tau),u(\xi(\tau),\tau),\dot{\xi}(\tau))d\tau\\
&\ \ \ \ -\phi(\bar{\gm}(0))-
\int_0^{t_1}L(\bar{\gm}(\tau),u(\bar{\gm}(\tau),\tau),\dot{\bar{\gm}}(\tau))d\tau,
\end{align*}which together with (\ref{xi}) gives rise to
\begin{equation}
u(\gm(t_2),t_2)-u(\gm(t_1),t_1)\leq \int_{t_1}^{t_2}L(\gm(\tau),u(\gm(\tau),\tau),\dot{\gm}(\tau))d\tau,
\end{equation}which verifies (i) of Definition \ref{nw}. By means of Lemma \ref{chc}, there exists a $C^1$ calibrated curve $\gm:[t_1,t_2]\rightarrow M$ with $\gm(t_2)=x$ such that
  \begin{equation}
u(x,t_2)-u(\gm(t_1),t_1)=\int_{t_1}^{t_2}L(\gm(\tau),u(\gm(\tau),\tau),\dot{\gm}(\tau))d\tau.
\end{equation}which implies (ii) of Definition \ref{nw}. This completes the proof of Lemma \ref{uisw}.
\End

Based on Definition \ref{visco}, it is easy to see that a variational solution of (\ref{hje}) is a viscosity solution.
\begin{Lemma}\label{awisv}
A variational solution of (\ref{hje}) with initial condition is a viscosity solution.
\end{Lemma}
\Proof
 Let $u$ be a variational solution of (\ref{hje}). Since $u(x,0)=\phi(x)$ it suffices to consider $t\in (0,T]$.  We use $V$ to denote an open subset of $M$. Let $\phi:V\times[0,T]\rightarrow \R$ be a $C^1$ test function such that $u-\phi$ has a maximum at $(x_0,t_0)$. This means $\phi(x_0,t_0)-\phi(x,t)\leq u(x_0,t_0)-u(x,t)$. Fix $v\in T_{x_0}M$ and for a given $\delta>0$, we choose a $C^1$ curve $\gm:[t_0-\delta,t_0+\delta]\rightarrow M$ with $\gm(t_0)=x_0$ and $\dot{\gm}(t_0)=\xi$. For $t\in [t_0-\delta,t_0]$, we have
\begin{align*}
\phi(\gm(t_0),t_0)-\phi(\gm(t),t)&\leq u(\gm(t_0),t_0)-u(\gm(t),t),\\
&\leq \int_t^{t_0}L(\gm(\tau),u(\gm(\tau),\tau),\dot{\gm}(\tau))d\tau,
\end{align*}where the second inequality is based on (i) of Definition \ref{weakkam}. Hence,
\begin{equation}
\frac{\phi(\gm(t),t)-\phi(\gm(t_0),t_0)}{t-t_0}\leq\frac{1}{t-t_0}
\int_{t_0}^tL(\gm(\tau),u(\gm(\tau),\tau),\dot{\gm}(\tau))d\tau.
\end{equation}Let $t\rightarrow t_0$, we have
\[\partial_t\phi(x_0,t_0)+\partial_{x}\phi(x_0,t_0)\cdot \xi\leq L(x_0,u(x_0,t_0),\xi),\]which together with Legendre transformation implies
\[\partial_t\phi(x_0,t_0)+H(x_0,u(x_0,t_0),\partial_x\phi(x_0,t_0))\leq 0,\]which shows that $u$ is a viscosity subsolution.

To complete the proof of Lemma \ref{awisv}, it remains to show that $u$ is a supersolution. $\psi:V\times[0,T]\rightarrow \R$ be a $C^1$ test function and $u-\psi$ has a minimum at $(x_0,t_0)$. We have $\psi(x_0,t_0)-\psi(x,t)\geq u(x_0,t_0)-u(x,t)$. From (ii) of Definition \ref{weakkam}, there exists a $C^1$ curve $\gm:[0,t_0]\rightarrow M$ with $\gm(t_0)=x_0$ and $\dot{\gm}(t_0)=\eta$ such that for $0\leq t<t_0$, we have
  \begin{equation}
u(\gm(t_0),t_0)-u(\gm(t),t)=\int_{t}^{t_0}L(\gm(\tau),u(\gm(\tau),\tau),\dot{\gm}(\tau))d\tau.
\end{equation}Hence
\[\psi(x_0,t_0)-\psi(x,t)\geq\int_{t}^{t_0}L(\gm(\tau),u(\gm(\tau),\tau),\dot{\gm}(\tau))d\tau.\]
Moreover, we have
\[\frac{\psi(\gm(t),t)-\psi(\gm(t_0),t_0)}{t-t_0}\geq\frac{1}{t-t_0}
\int_{t_0}^tL(\gm(\tau),u(\gm(\tau),\tau),\dot{\gm}(\tau))d\tau.\]Let $t$ tend to $t_0$, it gives rise to
\[\partial_t\psi(x_0,t_0)+\partial_{x}\psi(x_0,t_0)\cdot \eta\geq L(x_0,u(x_0,t_0),\eta),\]which implies
\[\partial_t\phi(x_0,t_0)+H(x_0,u(x_0,t_0),\partial_x\phi(x_0,t_0))\geq 0.\]This finishes the proof of Lemma \ref{awisv}.
\End

By the comparison theorem (see \cite{Ba3} for instance), it yields that  the viscosity solution of (\ref{hje}) is unique under the assumptions (H1)-(H5). So far, we have obtained that there exists a unique  viscosity solution $u(x,t)$ of (\ref{hje}) with initial condition $u(x,0)=\phi(x)$. Moreover, $u(x,t)$  can be represented implicitly as
\begin{equation}
u(x,t)=\inf_{\gm(t)=x}\left\{\phi(\gm(0))+\int_0^tL(\gm(\tau),u(\gm(\tau),\tau),\dot{\gm}(\tau))d\tau\right\}.
\end{equation}This completes the proof of Theorem \ref{two}.

\section{\sc Solution semigroup}
In this section,
we will prove  $T_t$ is a semigroup, which is called a solution semigroup. A similar definition was also introduced by  \cite{D} under more strict conditions on $H$. Under the assumptions (H1)-(H5), we will detect some further properties of the solution semigroup. Moreover, we will complete the proof of Theorem \ref{four}.

\subsection{Semigroup property of $T_t$}
Based on Theorem \ref{one},  we have
\begin{equation}\label{fixu1}
u(x,t)=\inf_{\gm(t)=x}\left\{\phi(\gm(0))+\int_0^tL(\gm(\tau),u(\gm(\tau),\tau),\dot{\gm}(\tau))d\tau\right\},
\end{equation}where the infimums are taken among absolutely continuous curves. In particular, the infimums are attained at the characteristics of (\ref{hje}). From (\ref{fixu1}), it follows that $u(x,t)$ can be represented as
\[u(x,t)=T_t\phi(x),\]
where $T_t$ denotes a nonlinear operator. Hence, we have
\begin{equation}\label{51}
T_t\phi(x)=\inf_{\gm(t)=x}\left\{\phi(\gm(0))+\int_0^tL(\gm(\tau),T_\tau\phi(\gm(\tau)),
\dot{\gm}(\tau))d\tau\right\}.
\end{equation}The following lemma implies $T_t$ is a semigroup.
\begin{Lemma}\label{semi}
$\{T_t\}_{t\geq 0}$ is a one-parameter semigroup of operators from
$C(M,\R)$ into itself.
\end{Lemma}
\Proof
It is easy to see $T_0  = Id$. It suffices to prove that
$T_{t+s}=T_t\circ T_s$ for any $t,s\geq 0$.

For every $\eta\in C(M, \mathbb{R})$ and $u \in C(M\times
[0,T],\mathbb{R})$, we define an operator $\A_t^{\eta}$ such that
\begin{equation}
\A^{\eta}[u](x,t) = \inf_{\gamma(t)=x} \left\{ \eta\big(\gamma(0)\big) + \int_0^t
L\big(\gamma(\tau),u\big(\gamma(\tau),\tau\big), \dot\gamma(\tau)\big) d\tau
\right\}.
\end{equation}
By virtue of Theorem \ref{one}, it follows that $\A^{\eta}$ has a unique fixed point.

By (\ref{51}), we have
\begin{equation*}
\begin{split}
T_t\circ T_s \phi(x) &= \inf_{\gamma(t)=x} \left\{
T_s\phi\big({\gamma}(0)\big) + \int_0^{t}
L\big({\gamma}(\tau),T_{\tau}\circ T_s
\phi\big({\gamma}(\tau)\big), \dot{{\gamma}}(\tau)\big)d\tau
\right\}\\
& = \A^{T_s \phi} [T_t\circ T_s \phi](x).
\end{split}
\end{equation*}

On the other hand,
\begin{equation*}
\begin{split}
T_{t+s} \phi(x) &= \inf_{\substack{\gamma(t+s)=x}} \left\{ \phi\big(\gamma(0)\big) +
\int_0^{t+s} L\big(\gamma(\tau),T_\tau
\phi\big(\gamma(\tau)\big),
\dot\gamma(\tau)\big)d\tau \right\}\\
& = \inf_{\substack{\gamma(t+s)=x}}
\left\{ \phi\big(\gamma(0)\big) + \bigg(\int_0^{s} +
 \int_{s}^{t+s}\bigg) L\big(\gamma(\tau),T_\tau
\phi\big(\gamma(\tau)\big), \dot\gamma(\tau)\big) d\tau \right\}\\
& = \inf_{\substack{\gamma(t+s)=x}}
\left\{ T_s\phi\big(\gamma(s)\big) + \int_s^{t+s}
L\big(\gamma(\tau),T_\tau \phi\big(\gamma(\tau)\big),
\dot\gamma(\tau)\big)d\tau \right\}\\
& = \inf_{\substack{\bar{\gamma}(t)=x}} \left\{ T_s\phi\big(\bar{\gamma}(0)\big) +
\int_0^{t} L\big(\bar{\gamma}(\tau),T_{\tau+s}
\phi\big(\bar{\gamma}(\tau)\big),
\dot{\bar{\gamma}}(\tau)\big)d\tau \right\}  \\
& = \A^{T_s \phi} [T_{t+s} \phi](x).
\end{split}
\end{equation*}
Hence, both $T_t\circ T_s \phi$ and $T_{t+s} \phi$ are fixed points of $\A^{T_s \phi} $, which together with the uniqueness of the fixed point of $\A^{T_s \phi} $ yields  $T_{t+s}=T_t\circ T_s$. This completes
the proof of  Lemma \ref{semi}.
\End
To fix the notion, we call $T_t$ a solution semigroup. In the following subsections, we will prove some further properties of the solution semigroup $T_t$.
\subsection{Properties of the solution semigroup}
First of all, it is easy to obtain the following proposition about the monotonicity of $T_t$.
\begin{Proposition}[Monotonicity]\label{monotonicity}
For given $\phi,\psi\in C(M,\R)$ and $t\geq 0$, if $\phi\leq\psi$, then $T_t\phi\leq T_t\psi$.
\end{Proposition}
\Proof
 For given $\phi,~\psi\in C(M,\mathbb{R})$ with
$\phi\leq \psi$, by contradiction, we assume that there exist
$t_1>0$ and $x_1\in M$ such that $T_{t_1}\phi(x_1) > T_{t_1}
\psi(x_1)$.
Let $\gm_\psi: [0,t_1]\rightarrow M$ be a calibrated curve of $T_t\psi$ with $\gm_\psi(t_1)=x_1$. We denote
\[F(\tau)=T_{\tau}\phi(\gm_\psi(\tau)) - T_{\tau}
\psi(\gm_\psi(\tau)).\] It is easy to see that $F(\tau)$ is continuous and $F(t_1)>0$. Since
 \[F(0)=T_{0}\phi(\gm_\psi(0) - T_{0}
\psi(\gm_\psi(0))\leq 0,\]
there exists $t_0\in [0,t_1)$ such that $F(t_0)=0$ and for any $\tau\in [t_0,t_1]$, $F(\tau)\geq 0$, i.e.
\begin{equation}\label{52}
T_{\tau}\phi(\gm_\psi(\tau)) \geq T_{\tau}
\psi(\gm_\psi(\tau)).
\end{equation}

Moreover, it follows from (\ref{51}) that
\begin{equation}\label{monotonicity estimate}
\begin{split}
&\quad T_{t_1}\phi(x_1) - T_{t_1}\psi(x_1) \\&= \inf_{\substack{\gamma(t_1)=x_1}} \left\{ T_{t_0}\phi\big(\gamma({t_0})\big)
+ \int_{t_0}^{t_1} L\big(\gamma(\tau),T_\tau
\phi\big(\gamma(\tau)\big),
\dot\gamma(\tau)\big) \ d\tau \right\} \\
&\quad-\inf_{\substack{\gamma(t_1)=x_1\\}}
\left\{ T_{t_0}\psi\big(\gamma({t_0})\big) + \int_{t_0}^{t_1}
L\big(\gamma(\tau),T_\tau
\psi\big(\gamma(\tau)\big), \dot\gamma(\tau)\big) \ d\tau \right\} \\
& \leq T_{t_0} \phi(\gm_\psi({t_0})) - T_{t_0} \psi(\gm_\psi({t_0})) \ +\\
&\quad \int_{t_0}^{t_1}\left( L\big(\gm_\psi(\tau),T_\tau
\phi\big(\gm_\psi(\tau)\big), \dot\gm_\psi(\tau)\big) -
L\big(\gm_\psi(\tau),T_\tau \psi\big(\gm_\psi(\tau)\big),
\dot{\gm}_\psi(\tau)\big)\right) \ d\tau,\\
&\leq\int_{t_0}^{t_1}\left( L\big(\gm_\psi(\tau),T_\tau
\phi\big(\gm_\psi(\tau)\big), \dot\gm_\psi(\tau)\big) -
L\big(\gm_\psi(\tau),T_\tau \psi\big(\gm_\psi(\tau)\big),
\dot{\gm}_\psi(\tau)\big)\right) \ d\tau.
\end{split}
\end{equation}Combining with (\ref{52}) and (L5), we have
\[
 T_{t_1}\phi(x_1) \leq T_{t_1}\psi(x_1),
\]
which is a contradiction. This finishes the proof of Proposition \ref{monotonicity}.\End

By a similar argument as the one in Proposition \ref{monotonicity}, one can obtain the non-expansiveness of $T_t$. For $\phi\in C(M,\R)$, we use $\|\phi\|_\infty$ to denote $C^0$-norm of $\phi$. We have the following proposition.
\begin{Proposition}[Non-expansiveness]\label{nonex}
For given $\phi,\psi\in C(M,\R)$ and $t\geq 0$, we have $\|T_t\phi-T_t\psi\|_\infty\leq \|\phi-\psi\|_\infty$.
\end{Proposition}
\Proof
By contradiction, we assume that there exist
$t_1>0$ and $x_1\in M$ such that $T_{t_1}\phi(x_1) -T_{t_1}
\psi(x_1)> \|\phi-\psi\|_\infty$.
Let $\gm_\psi: [0,t_1]\rightarrow M$ be a calibrated curve of $T_t\psi$ with $\gm_\psi(t_1)=x_1$. We denote
\[G(\tau)=T_{\tau}\phi(\gm_\psi(\tau)) - T_{\tau}
\psi(\gm_\psi(\tau))-\|\phi-\psi\|_\infty.\] It is easy to see that $G(\tau)$ is continuous and $G(t_1)>0$.
Since
 \[G(0)=T_{0}\phi(\gm_\psi(0) - T_{0}
\psi(\gm_\psi(0))-\|\phi-\psi\|_\infty\leq 0,\]
there exists $t_0\in [0,t_1)$ such that $G(t_0)=0$ and for any $\tau\in [t_0,t_1]$, $G(\tau)\geq 0$, i.e.
\begin{equation}\label{53}
T_{\tau}\phi(\gm_\psi(\tau)) - T_{\tau}
\psi(\gm_\psi(\tau))\geq\|\phi-\psi\|_\infty.
\end{equation}
A similar calculation as (\ref{monotonicity estimate}) implies
\begin{equation}\label{54}
T_{t_1}\phi(x_1) -T_{t_1}
\psi(x_1)\leq\|\phi-\psi\|_\infty
\end{equation}which is a contradiction. Hence, we have
\[T_{t_1}\phi(x_1) -T_{t_1}
\psi(x_1)\leq\|\phi-\psi\|_\infty.\]

Similarly, we have
\[T_{t_1}\phi(x_1) -T_{t_1}
\psi(x_1)> -\|\phi-\psi\|_\infty.\]
 This finishes the proof of Proposition \ref{monotonicity}.
\End

We use $c(L(x,a,\dot{x}))$ to denote the Ma\~{n}\'{e} critical value of $L(x,a,\dot{x})$. By \cite{CIPP}, we have
\begin{equation}
c(L(x,a,\dot{x}))=\inf_{u\in C^1(M,\R)}\sup_{x\in M}H(x,a,\partial_xu).
\end{equation}
Let
\begin{equation}
\mathcal{C}=\left\{c(L(x,a,\dot{x})): a\in\R\right\}.
\end{equation}It follows that for any $c\in \mathcal{C}$, there exists $a\in\R$ such that $c(L(x,a,\dot{x}))=c$. Let
$L_c=L+c$, then $c(L_c(x,a,\dot{x}))=0$. In the following context, we consider $L_c$ instead of $L$. Moreover,
\begin{equation}
T_t\phi(x)=\inf_{\gm(t)=x}\left\{\phi(\gm(0))+\int_0^tL_c(\gm(\tau),T_\tau\phi(\gm(\tau)),
\dot{\gm}(\tau))d\tau\right\}.
\end{equation}

 Without ambiguity, we still use $L$ to denote $L_c$ for the simplicity of notations. The following proposition implies $T_t\phi$ is uniformly bounded.
\begin{Proposition}[Uniform bound]\label{uniform bounds for semigroup}
For every $\phi\in C(M,\mathbb{R})$, there exists
a positive constant $K$ such that for any $t\geq 0$
\begin{equation}\label{uk}
\|T_t \phi\|_\infty \leq K.
\end{equation}
\end{Proposition}
\Proof
For $t=0$, $T_0\phi=\phi$, which is uniformly bounded. Let $u(x,t):=T_t\phi(x)$.

On the one hand, we show that $u(x,t)$ is uniformly bounded from below. Without loss of generality, one can assume $u(x,t) < \a$. For every $(x,t)\in M\times[0,+\infty)$, there exists a calibrated curve $\gm$ of $u$ with $\gm(t)=x$ such that
\begin{equation}
u(x,t)=\phi(\gm(0))+\int_0^tL(\gm(\tau),u(\gm(\tau),\tau),
\dot{\gm}(\tau))d\tau.
\end{equation}
Then, we have the following two cases:
\begin{itemize}
\item [(I)] there exists a $\tau_0\in [0,t)$ such that $u(\gm(\tau_0),\tau_0) = \a$ and $u(\gm(\tau),\tau)< \a$ for $\tau\in [\tau_0,t]$;
\item [(II)] for every $\tau\in[0,t]$,  $u(\gm(\tau),\tau) < \a$.
\end{itemize}
 For Case (I), it follows from (L5) that
  \begin{equation*}
  \begin{split}
  u(x, t) & =u(\gm(\tau_0),\tau_0) + \int_{\tau_0}^t L\big( \gamma(\tau), u(\gm(\tau),\tau), \dot{\gamma}(\tau)\big)\ d\tau,\\
  &\geq \a + \int_{\tau_0}^t L\big( \gamma(\tau), \a, \dot{\gamma}(\tau)\big)\ d\tau\\
  &\geq \a +  h_{t-\tau_0}(\gm(\tau_0),x).
  \end{split}
  \end{equation*}where $h_{t-\tau_0}(\gm(\tau_0),x)$ denotes the minimal action from $\gm(\tau_0)$ to $x$ for the Lagrangian $L(x, \a, \dot{x})$. It is easy to see that $h_{t-\tau_0}(\gm(\tau_0),x)$ is uniformly bounded from below. Hence,  $u(x, t)$ is uniformly bounded from below.

 For Case (II), a similar calculation yields
 \[u(x,t)\geq \min_{x\in M} \phi(x) +  h_t (x,y).\]

It follows from the compactness of $M$ that there exists a constat $K_1$ independent of $(x,t)$ such that  $u(x,t)\geq K_1$ for any $(x,t) \in M\times[0,+\infty)$.

On the other hand, we show that $u(x,t)$ is uniformly bounded from above. Without loss of generality, one can assume $u(x,t) > \a$.  For a given $x_0\in M$, let $\bar{\gm}: [0,t]\rightarrow M$ be a minimal curve with $\bar{\gm}(0)=x_0$ and $\bar{\gm}(t)=x$ such that
\[\int_0^tL(\bar{\gm}(\tau),a,\dot{\bar{\gm}}(\tau))d\tau=\inf_{\substack{\gm(0)=x_0\\ \gm(t)=x}}\int_0^tL(\gm(\tau),a,\dot{\gm}(\tau))d\tau,\]
where the infimums are taken among the absolutely continuous curves.
Then, we have the following two cases:
\begin{itemize}
\item [(I)] there exists a $\tau_0\in [0,t)$ such that $u(\bar{\gm}(\tau_0),\tau_0) = \a$ and $u(\bar{\gm}(\tau),\tau)> \a$ for $\tau\in [\tau_0,t]$;
\item [(II)] for every $\tau\in[0,t]$,  $u(\bar{\gm}(\tau),\tau) > \a$.
\end{itemize}
 For Case (I), we have
\begin{equation*}
\begin{split}
u(x,t) &\leq u(\bar{\gm}(\tau_0), \tau_0) + \int_{\tau_0}^t L\big( \bar{\gm}(\tau), u(\bar{\gm}(\tau),\tau), \dot{\bar{\gm}}(\tau)\big)\ d\tau\\
       &\leq \a + \int_{\tau_0}^t L\big( \bar{\gm}(\tau), \a, \dot{\bar{\gm}}(\tau)\big)\ d\tau\\
       &= \a + h_{t-\tau_0}(\bar{\gm}(\tau_0), x).
\end{split}
\end{equation*}
To verify that $u(x,t)$ is bounded above, it suffices to prove that $h_{t-\tau_0}( \bar{\gm}(\tau_0), x)$ is bounded from above.
Without loss of generality, we assume $t>1$. Hence, among $\tau_0$ and $t-\tau_0$, there exists at least one not less than $\frac{1}{2}$. If $t-\tau_0\geq \frac{1}{2}$, then it follows  from (\ref{lin}) that $h_{t-\tau_0}(\bar{\gm}(\tau_0),x)$ is bounded from above. If $\tau_0\geq \frac{1}{2}$, it follows from
 (\ref{wang}) that
\begin{equation*}
h_{t- \tau_0}(\bar{\gm}(\tau_0), x)=h_{t}(x_0,x)-h_{\tau_0}(x_0, \bar{\gm}(\tau_0)).
\end{equation*}which also implies $h_{t-\tau_0}(\bar{\gm}(\tau_0),x)$ is uniformly bounded from above. Hence,  $u(x, t)$ is uniformly bounded from above.

For Case (II), by a similar argument, we can obtain the upper bound of $h_{t}(x,y)$. Therefore, there exists a constant $K_2$ independent of $(x,t)$ such that $u(x,t)\leq K_2$ for any $(x,t) \in M\times[0,+\infty)$.
This completes the proof of Proposition \ref{uniform bounds for semigroup}.
\End

Based on Proposition \ref{uniform bounds for semigroup}, we can obtain the equi-Lipschitz of the familiy of functions $T_t\phi(x)$.
\begin{Proposition}[Equi-Lipschitz]\label{equili}
For every $\phi(x)\in C(M,\R)$ and $\delta>0$, the family of functions $(x,t)\rightarrow T_t\phi(x)$ is equi-Lipschitz on $(x,t)\in M\times [\delta,+\infty)$.
\end{Proposition}
The key point to prove Proposition \ref{equili} is a priori compactness from which it is easy to verify Proposition \ref{equili} following from a similar argument as \cite{F3}. Let $u(x,t):=T_t\phi(x)$. We are devoted to proving the following lemma.
\begin{Lemma}[A Priori Compactness]\label{apc}
For a given $\delta>0$, there exists a compact subset $\mathcal{K}_\delta$ such that for every calibrated curve $\gm$ of $u$ and any $t>\delta$, we have
\begin{equation*}
(\gm(t),u(\gm(t),t),\dot{\gm}(t))\in \mathcal{K}_\delta.
\end{equation*}
\end{Lemma}
\Proof By Proposition \ref{uniform bounds for semigroup} and the compactness of $M$, it suffices to prove that there exists  a constant $A>0$ such that for any $t>\delta$
\begin{equation*}
|\dot{\gamma}(t)| \leq A.
\end{equation*}

We assume  by contradiction that for any $n\in\N$, there exists a $t_n>\delta$  such that
the calibrated curve
\begin{equation}\label{priori}
|\dot{\gamma}(t_n)|\geq n.
\end{equation}
Based on (\ref{priori}), we have the following claim.

\textbf{Claim:} For any $m>0$, if $n$ is large enough, then for any $\tau\in [t_n-\delta,t_n]$, we have
\begin{equation}\label{priori1}
|\dot{\gamma}(\tau)|\geq m.
\end{equation}

\Proof
We assume  by contradiction that there exists $m_0>0$ such that for any $n\in\N$, one can find a sequence $\tau_n\in [t_n-\delta,t_n]$ satisfying $|\dot{\gamma}(\tau_n)|\leq m_0$. Let $s_n:=\tau_n-(t_n-\delta)$, then $s_n\in [0,\delta]$. Extracting a subsequence if necessary, it follows from (\ref{uk}), (\ref{priori1}) and the compactness of $M$ that
\[(\gm(\tau_n),u(\gm(\tau_n),\tau_n),\dot{\gm}(\tau_n))\rightarrow (\bar{x},\bar{u},\bar{v}),\quad s_n\rightarrow \bar{s}.\]
Let $\Phi_s$ be the flow generated by $L(x,u,\dot{x})$. Then it follows from the completeness of the flow that for any $s\in [0,\delta]$, $\Phi_{s}(\bar{x},\bar{u},\bar{v})$ is well defined. Theorem \ref{one} implies $(\gm(s),u(\gm(s),s),\dot{\gm}(s))$ is the flow generated by $L(x,u,\dot{x})$. We consider $(\bar{x},\bar{u},\bar{v})$ as the initial condition of $\Phi_s$. Based on the continuous dependence of solutions of ODEs on initial conditions, it follows that for $n$ large enough,
\[|(\gm(t_n),u(\gm(t_n),t_n),\dot{\gm}(t_n))-\Phi_{\delta-\bar{s}}(\bar{x},\bar{u},\bar{v})|\leq \epsilon,\]
where $\epsilon$ is a small constant independent of $n$.
By virtue of (L3), it yields that $|\dot{\gm}(t_n)|$ has a bound independent of $n$, which is in contradiction with (\ref{priori}). This verifies the claim.
\End

By Proposition \ref{uniform bounds for semigroup}, we have
$|u(x,t)|\leq K$. It follows from the superlinear growth (see (L2)) that there exists $C>0$ such that
\begin{equation}\label{62}
L(x, K, \dot{x})\geq |\dot{x}|-C.
\end{equation}
From (\ref{priori1}) and (\ref{62}), it follows that
\begin{equation*}
\begin{split}
u(\gm &(t_n),t_n) - u(\gamma(t_n-\delta), t_n-\delta)= \int_{t_n-\delta}^{t_n} L(\gamma(\tau), u(\gm(\tau),\tau), \dot{\gamma}(\tau))\ d\tau\\
&\geq \int_{t_n-\delta}^{t_n} L(\gamma(\tau), K, \dot{\gamma}(\tau))\ d\tau\geq \int_{t_n-\delta}^{t_n}|\dot{\gm}(\tau)|\ d\tau-C\delta,\\
                          &\geq \frac{\delta}{2}m-C\delta,
\end{split}
\end{equation*}where the second inequality is owing to the assumption (L5). On the other hand, we have
\begin{equation}
|u(\gm (t_n),t_n) - u(\gamma(t_n-\delta), t_n-\delta)|\leq 2K.
\end{equation}
Hence, we have
\[\frac{\delta}{2}m-C\delta\leq 2K,\]
which is a contradiction for $m$ large enough. In fact, it suffices to take $m>\frac{4K}{\delta}+2C$. Hence, there exists  a constant $A>0$ such that for any $t>\delta$,
$|\dot{\gamma}(t)| \leq A$. This completes the proof of Lemma \ref{apc}.
\End
Based on Lemma \ref{apc}, Proposition \ref{equili} can be obtained following from a similar lengthy and tedious argument as \cite{F3}, where we omit it.

So far, we complete the proof of Theorem \ref{four}.

\section{\sc Convergence of the solution semigroup}
In this section, we are concerned with the convergence of the solution semigroup generated by the Lagrangian $L_c:=L+c$.  We still use $L$ instead of $L_c$ in the following. We will show  that for any $\phi(x)\in C(M,\R)$, $T_t\phi$
converges as $t\rightarrow \infty$ to a weak KAM solution of
\begin{equation}\label{stahh}
H(x, u, \partial_x u)=0.
\end{equation}
Moreover, we will complete the proof of Theorem \ref{five}.

\subsection{Step 1}
In this step, we will prove the existence of a weak KAM solution of (\ref{stahh}). By virtue of Proposition \ref{uniform bounds for semigroup}, we have $T_t\phi$ is uniformly bounded for any $\phi\in C(M,\R)$. Hence, $\limsup_{t\rightarrow\infty} T_t\phi$ does exist, which is denoted by $\bar{u}$. We have the following lemma.

\begin{Lemma}\label{common fixed point existence}
The limit $\lim_{t\rightarrow\infty}T_t\bar{u}$ does exist. Moreover,
let
\[u_\infty(x):=\lim_{t\rightarrow\infty}T_t\bar{u}(x),\]
then $u_\infty(x)$ is a weak KAM solution of (\ref{stahh}).
\end{Lemma}
\Proof
Due to the definition of limsup, for every
$\epsilon>0$, there exists  $s_0\in\mathbb{R}^+$ such that for any $s\geq s_0$, we have
\begin{equation}
T_s\phi\leq \bar{u}+\epsilon,
\end{equation}
which the non-expansiveness and monotonicity of $T_t$ implies
\begin{equation}
T_t\circ T_s \phi \leq  T_t(\bar{u}+\epsilon)\leq T_t\bar{u}+\epsilon.
\end{equation}
Fixing $t\geq 0$, we take limsup for the above inequality as $s\rightarrow\infty$. Since
\begin{equation}
\limsup_{s\rightarrow\infty}T_t\circ T_s\phi=\limsup_{t+s\rightarrow\infty}T_{t+s}\phi=\bar{u},
\end{equation}
then we
obtain
\begin{equation}
\bar{u}\leq T_t \bar{u}+\epsilon.
\end{equation}
Since $\epsilon$ is arbitrary, we have
\[
\bar{u}\leq T_t \bar{u}.\]
 Hence, by the monotonicity of $T_t$, it follows from the semigroup property that $T_t \bar{u}$ is non-decreasing with respect to $t$. Combining with boundedness of $T_t\bar{u}$, it follows that the limit $\lim_{t\rightarrow\infty}T_t\bar{u}$ does exist, which is denoted by $u_\infty$. Then, we have
 \begin{equation}
 T_tu_\infty=u_\infty.
 \end{equation}
Based on Proposition \ref{uniform bounds for semigroup} and Proposition \ref{equili}, it follows from
Arzela-Ascoli theorem that
$u_\infty(x)\in C(M,\R)$.

It remains to verify $u$ is a weak KAM solution of (\ref{stahh}). By virtue of a similar argument as Lemma \ref{uisw},  it yields that for each continuous piecewise $C^1$ curve $\gm:[t_1,t_2]\rightarrow M$ where $0\leq t_1<t_2\leq T$, we have
\begin{equation}
u_\infty(\gm(t_2))-u_\infty(\gm(t_1))\leq\int_{t_1}^{t_2}L(\gm(\tau),
u_\infty(\gm(\tau)),\dot{\gm}(\tau))d\tau,
\end{equation}which implies (i) of Definition \ref{weakkam}. In addition,
there exists a $C^1$ calibrated curve $\gm_t:[-t,0]\rightarrow M$ with $\gm_t(0)=x$ such that for any $t'\in [-t,0]$, we have
  \begin{equation}
u_\infty(x)-u_\infty(\gm_t(t'))=\int_{t'}^{0}L(\gm_t(\tau),u_\infty(\gm_t(\tau)),\dot{\gm}_t(\tau))d\tau.
\end{equation}
Based on the a priori compactness given by Lemma \ref{apc}, for a given $\delta>0$, there exists a compact subset $\mathcal{K}_\delta$ such that for  any $s>\delta$, we have
\begin{equation*}
(\gm_t(s),u_\infty(\gm_t(s)),\dot{\gm}_t(s))\in \mathcal{K}_\delta.
\end{equation*}
Since $\gm_t$ is a calibrated curve, it follows from Lemma \ref{chc} that
\[(\gm_t(s),u_\infty(\gm_t(s)),\dot{\gm}_t(s))=\Phi_s(\gm_t(0),u_\infty(\gm_t(0)),\dot{\gm}_t(0))=\Phi_s(x,u_\infty(x),\dot{\gm}_t(0)).\]
The points $(\gm_t(0),u_\infty(\gm_t(0)),\dot{\gm}_t(0))$ are contained in a compact subset, then one can find a sequence $t_n$ such that $(x,\dot{\gm}_{t_n}(0))$ tends to $(x,v_\infty)$ as $n\rightarrow \infty$. Fixing $t'\in (-\infty,0]$, the function $s\mapsto \Phi_s(x,u_\infty(x),\dot{\gm}_{t_n}(0))$ is defined on $[t',0]$ for $n$ large enough. By the continuity of $\Phi_s$, the sequence converges uniformly on the compact interval $[t',0]$ to the map $s\mapsto \Phi_s(x,v_\infty)$. Let
\[(\gm_\infty(s),u_\infty(\gm_\infty(s)), \dot{\gm}_\infty(s)):=\Phi_s(x,v_\infty),\] then
for any $t'\in (-\infty,0]$, we have
  \begin{equation}
u_\infty(x)-u_\infty(\gm_\infty(t'))=\int_{t'}^{0}L(\gm_\infty(\tau),u_\infty(\gm_\infty(\tau)),\dot{\gm}_\infty(\tau))d\tau,
\end{equation}
which implies (ii) of Definition \ref{weakkam}.
Hence, $u_\infty$ is a weak KAM solution of (\ref{stahh}). This completes the proof of Lemma \ref{common fixed point existence}.
\End

\subsection{Step 2}
Since $u_\infty(x)$ is a weak KAM solution, then $u_\infty(x)$ is Lipschitz. Let $\mathcal{D}$ be the set of all differentiable points of $u_\infty$ on $M$. Due to the Lipschitz property of $u_\infty$, it follows that $\mathcal{D}$ has full Lebesgue measure. For $x\in \mathcal{D}$, we have
\begin{equation}\label{hxud}
H(x, u_\infty(x),\partial_x u_\infty(x) ) = 0.
\end{equation}
We define
\begin{equation}\label{define}
\widetilde{L}(x,\dot{x})= L(x, u_\infty(x),\dot{x}) - \langle \partial_xu_\infty(x), \dot{x}\rangle,\quad  x\in \mathcal{D}.
\end{equation}
Denote
\begin{equation}\label{gamm}
\Gamma :=\left\{ \left(x, \frac{\partial H}{\partial p}(x,u_\infty(x),\partial_x u_\infty(x))\right)~:~x\in \mathcal{D}  \right\},
\end{equation}
where $\frac{\partial H}{\partial p}$ denotes the partial derivative of $H$ with respect to the third argument. Hence, we have the following lemma.
\begin{Lemma}\label{squr}
For any $x\in\mathcal{D}$, $\widetilde{L}(x,\dot{x})\geq 0$. In particular, $\widetilde{L}(x,\dot{x})=0$ if and only if $(x,\dot{x})\in \Gamma$.
\end{Lemma}
\Proof
By (\ref{define}) and (\ref{gamm}), we have
\begin{equation}\label{11}
\widetilde{L}\big|_\Gamma =-H(x, u_\infty(x),\partial_x u_\infty(x))= 0.
\end{equation}
In addition, we have
\begin{equation}\label{12}
\frac{\partial \widetilde{L}}{\partial \dot{x}}\big|_\Gamma
= \frac{\partial L}{\partial \dot{x}}(x,u_\infty(x),\dot{x}) - \partial_x u_\infty(x) = 0.
\end{equation}

By (L2), it follows  from (\ref{11}) that there exists $K>0$ large enough such that for $|\dot{x}|> K$,
\[\widetilde{L} (x, \dot{x})\geq d>0,\]where $d$ is a constant independent of $(x,\dot{x})$.

For $x\in \mathcal{D}$, $u_\infty(x)$ satisfies the equation (\ref{hxud}). It follows from (L2) that  $\partial_xu_\infty(x)$ is bounded. Let
\[\dot{x}_0:=\frac{\partial H}{\partial p}(x,u_\infty(x),\partial_x u_\infty(x)),\]
 then there exists $K'>0$ independent of $x$ such that $|\dot{x}_0|\leq K'$. Without loss of generality, one can assume $K'<K$.  From the assumption (L1), it follows that
$\frac{\partial^2L}{\partial \dot{x}^2}(x,u_\infty(x),\dot{x})$ is positive definite. Hence, for $|\dot{x}|\leq K$, it follows  from (\ref{11}) and (\ref{12}) that there exists $\Lambda>0$ independent of $(x,\dot{x})$ such that
\begin{equation}
\widetilde{L} (x, \dot{x}) \geq \Lambda\left|\dot{x} - \frac{\partial H}{\partial p}(x,u_\infty(x),\partial_x u_\infty(x))\right|^2.
\end{equation}
Consequently, it is easy to see that
\begin{equation}\label{Lagrangian graph}
\widetilde{L}(x, \dot{x})
\left\{\!\!\!
  \begin{array}{rl}
   &=0,  \qquad (x,\dot{x})\in \Gamma,\\
   &>0,  \qquad (x,\dot{x}) \notin \Gamma.
  \end{array}
\right.
\end{equation}This completes the proof of Lemma \ref{squr}.
\End
\subsection{Step 3}
In this step, we focus on the evolution of $H$ along the characteristics. First of all, we recall the definition of $T_t$ :
\begin{equation}\label{71}
T_t\phi(x)=\inf_{\gm(t)=x}\left\{\phi(\gm(0))+\int_0^tL(\gm(\tau),T_\tau\phi(\gm(\tau)),
\dot{\gm}(\tau))d\tau\right\},
\end{equation}where the infimums are taken among the absolutely continuous curves. In particular, the infimums are attained at  the characteristics of (\ref{hje}) based on Theorem \ref{one}. More precisely, let $\gm(s):[0,t]\rightarrow M$ be a calibrated curve of $T_t\phi$, then $(\gm(s),u(s),p(s))$ defined as
\begin{equation}
\left(\gm(s),u(s):=T_s\phi(\gm(s)),p(s)=\frac{\partial L}{\partial \dot{x}}(\gm(s),T_s\phi(\gm(s)),\dot{\gm}(s))\right)
\end{equation}
is $C^1$ and satisfies the characteristic equation (\ref{hje}). To avoid the ambiguity, we denote the characteristics by $(X(t),U(t),P(t))$.
Along the characteristics, a simple calculation yields that for any $s\in [0,t]$,
\begin{equation}\label{key}
\begin{split}
\frac{dH}{ds}(X(s),U(s),P(s)) &= \frac{\partial H}{\partial x}\
\dot{X}(s) + \frac{\partial H}{\partial u}\ \dot{U}(s) +
\frac{\partial H}{\partial p}\ \dot{P}(s)\\
& = -\frac{\partial H}{\partial u}(X(s),U(s),P(s))\
H(X(s),U(s),P(s)).
\end{split}
\end{equation}

Let $\bar{H}(s):= H(X(s),U(s),P(s))$. It follows from (H5) that
\begin{itemize}
\item  $\bar{H}(s)$ is a decreasing function with respect to $s$ if $\bar{H}(0) >0$;
\item  $\bar{H}(s)$ is an increasing function with respect to $s$ if $\bar{H}(0) < 0$;
\item  $\bar{H}(s)=0$ if $\bar{H}(0) = 0$.
\end{itemize}

By virtue of Proposition \ref{uniform bounds for semigroup} and Proposition \ref{equili}, it follows from Arzela-Ascoli theorem that
there exist a  sequence $t_n\rightarrow \infty$  and a Lipschitz function $\tilde{u}$
such that $T_{t_n} \phi \rightarrow \tilde{u}$. Based on Lemma \ref{common fixed point existence}, we have
\begin{equation}\label{m66}
\tilde{u}\leq u_\infty,
\end{equation}where $u_\infty$ denotes the weak KAM solution in Lemma \ref{common fixed point existence}.

For a given $s> 0$, extracting a subsequence if necessary, there exists a Lipschitz function $u^-(x)$ such that
\begin{equation}
\lim_{n\rightarrow\infty}T_{t_n-s}\phi(x)=u^-(x).
\end{equation}
In addition, we have $u^-\leq u_\infty$.
Let $w(x,t):=T_tu^-(x)$. By the monotonicity of $T_t$, it yields
\begin{equation}\label{33}
w(x,t)\leq u_\infty(x).
\end{equation}
 It follows from Theorem \ref{two} that $w(x,t)$ is a viscosity solution of the following equation:
\begin{equation}\label{55}
\begin{cases}
\partial_tw(x,t)+H(x,w(x,t),\partial_xw(x,t))=0,\\
w(x,0)=u^-(x).
\end{cases}
\end{equation}
In particular, we have $w(x,s)=\tilde{u}(x)$, where $\tilde{u}$ is the same as the one in (\ref{m66}).
It is easy to see  $\tilde{u}(x)$ is a Lipschitz function defined on $M$ (see \cite{F3}). Let $\mathcal{D}'$ be the set of all differentiable points of $\tilde{u}$ on $M$. Due to the Lipschitz property of $\tilde{u}$, it follows that $\mathcal{D}'$ has full Lebesgue measure. Let  $\mathcal{D}'':=\mathcal{D}\cap\mathcal{D}'$, then  $\mathcal{D}''$ also has full Lebesgue measure. For the simplicity of notations, we still use $\mathcal{D}$ instead of $\mathcal{D}''$. Hence, we have the following lemma.
\begin{Lemma}\label{uu}
For $x\in \mathcal{D}$, we have
\begin{equation}\label{91}
H(x,\tilde{u}(x),\partial_x\tilde{u}(x))= 0.
\end{equation}
\end{Lemma}
\Proof
By contradiction, we assume that there exists $x_0\in \mathcal{D}$  such that
\begin{equation}
H(x_0,\tilde{u}(x_0),\partial_x\tilde{u}(x_0))\neq 0.
\end{equation}
In the following, we consider only the case when
\begin{equation}\label{13}
H(x_0,\tilde{u}(x_0),\partial_x\tilde{u}(x_0))=\delta>0.
\end{equation}
since the case when $H(x_0,\tilde{u}(x_0),\partial_x\tilde{u}(x_0))=-\delta<0$ is similar.
From the Legendre transformation,
\begin{align*}
H(x_0,\tilde{u}(x_0),\partial_x\tilde{u}(x_0))&=\sup_{\dot{x}_0}\{\langle \dot{x}_0,\partial_x\tilde{u}(x_0)\rangle-L(x_0,\tilde{u}(x_0),\dot{x}_0)\},\\
&=\left\langle \frac{\partial L}{\partial \dot{x}}(x_0,\tilde{u}(x_0),\dot{x}_0),\dot{x}_0\right\rangle-L(x_0,\tilde{u}(x_0),\dot{x}_0),
\end{align*}
where
\begin{equation}
\partial_x\tilde{u}(x_0)=\frac{\partial L}{\partial \dot{x}}(x_0,\tilde{u}(x_0),\dot{x}_0).
\end{equation}
By (\ref{13}), it yields that
\begin{equation}\label{15}
L(x_0,\tilde{u}(x_0),\dot{x}_0)-\left\langle \frac{\partial L}{\partial \dot{x}}(x_0,\tilde{u}(x_0),\dot{x}_0),\dot{x}_0\right\rangle=-H(x_0,\tilde{u}(x_0),\partial_x\tilde{u}(x_0))=- \delta.
\end{equation}
We denote
\[\check{L}(x,w(x,t),\dot{x}):=L(x,w(x,t),\dot{x})-\left\langle \partial_xu_\infty(x),\dot{x}\right\rangle,\]
where $w(x,t)$ is a viscosity solution of (\ref{55}). In particular, $w(x_0,s)=\tilde{u}(x_0)$.
In terms of (\ref{define}), we have $\check{L}(x,u_\infty(x),\dot{x})=\widetilde{L}(x,\dot{x})$.
By (\ref{33}), we have $w(x,t)\leq u_\infty(x)$. Then it follows from (L5) that
\begin{equation}\label{666}
\check{L}(x,w(x,t),\dot{x})\geq \widetilde{L}(x,\dot{x}).
\end{equation}
By Lemma \ref{squr}, we have
\begin{equation}\label{geq0}
\check{L}(x,w(x,t),\dot{x})\geq 0.
\end{equation}

We denote
\begin{equation}
\hat{L}(x,w(x,t),\dot{x}):=L(x,w(x,t),\dot{x})-\left\langle \frac{\partial L}{\partial \dot{x}}(x,w(x,t),\dot{x}),\dot{x}\right\rangle.
\end{equation}
Let
\[\Delta:=\{(x,w(x,t),\dot{x}) : \hat{L}(x,w(x,t),\dot{x})=0,\ (x,t)\in M\times[0,s],\ \dot{x}\in T_xM\}.\]

\textbf{Claim:}  $\Delta$ is compact.

\Proof
By the Legendre transformation, we have
\[\hat{L}(x,w(x,t),\dot{x})=-H(x,w(x,t),p),\]where $\dot{x}=\frac{\partial H}{\partial p}(x,w(x,t),p)$.
If $(x,w(x,t),\dot{x})\in \Delta$, then $H(x,w(x,t),p)=0$. According to Proposition \ref{uniform bounds for semigroup}, $w(x,t)$ is uniformly bounded, which together with (H2) yields there exists $C$ independent of $(x,t)$ such that $|p|\leq C$. Moreover, $\dot{x}$ is contained in a compact set. Hence, the claim follows from the compactness of $M$.
\End
%
Let
\[\Sigma(\delta):=\{(x,w(x,t),\dot{x}) : \hat{L}(x,w(x,t),\dot{x})\leq-\delta,\ (x,t)\in M\times[0,s],\ \dot{x}\in T_xM\}.\]
Then we have the following claim.

\textbf{Claim:} For $(x,w(x,t),\dot{x})\in \Sigma(\delta)$,   there exists $\theta(\delta)>0$ such that \[\text{dist}((x,w(x,t),\dot{x}),\Delta)\geq \theta(\delta),\] where ``dist" denotes a distance induced by the Riemannian metric on $TM\times\R$.

\Proof
 By contradiction, we assume that for any $\epsilon>0$,
  \[\text{dist}((x,w(x,t),\dot{x}),\Delta)< \epsilon.\]
  Hence, there exists a sequence $(x_n,w(x_n,t_n),\dot{x}_n)$ contained in $\Sigma(\delta)$ such that for $n$ large enough,
    \[\text{dist}((x_n,w(x_n,t_n),\dot{x}_n),\Delta)< \frac{1}{n}.\]
    Extracting a subsequence if necessary, one obtain $(x_n,w(x_n,t_n),\dot{x}_n)\rightarrow(\bar{x},w(\bar{x},\bar{t}),\dot{\bar{x}})$.
   Based on the compactness of $\Delta$, we have  $(\bar{x},w(\bar{x},\bar{t}),\dot{\bar{x}})\in \Delta$. Then, it follows from the definition of $\Delta$ that
 \begin{equation}\label{16}
\hat{L}(\bar{x},w(\bar{x},\bar{t}),\dot{\bar{x}})=0.
\end{equation}
Since $\hat{L}(x_n,w(x_n,t_n),\dot{x}_n)\leq-\delta<0$, then  it follows from the continuity of $\hat{L}$ with respect to $x$ that
\begin{equation}
\hat{L}(\bar{x},w(\bar{x},\bar{t}),\dot{\bar{x}})\leq- \delta<0,
\end{equation}which is in contradiction with (\ref{16}). Hence, for $(x,w(x,t),\dot{x})\in \Sigma(\delta)$, there exists $\theta>0$ such that \[\text{dist}((x,w(x,t),\dot{x}),\Delta)\geq \theta(\delta),\] which verifies the claim.
\End

\textbf{Claim:} For $x\in \mathcal{D}$ and $(x,w(x,t),\dot{x})\in \Sigma(\delta)$, there exists $\delta'>0$ such that
\begin{equation}
\check{L}(x,w(x,t),\dot{x})\geq \delta'>0,
\end{equation}

\Proof
 By contradiction, we assume that for any $\epsilon>0$,
  \[\check{L}(x,w(x,t),\dot{x})< \epsilon.\]
  Hence, there exists a sequence $(x_n,w(x_n,t_n),\dot{x}_n)$ satisfying $x_n\in \mathcal{D}$ and $(x_n,w(x_n,t_n),\dot{x}_n)\in \Sigma$ such that for $n$ large enough,
    \[\check{L}(x_n,w(x_n,t_n),\dot{x}_n)< \frac{1}{n}.\] By the definition of $\check{L}$, we have
    \begin{equation}\label{555}
L(x_n,w(x_n,t_n),\dot{x}_n)-\left\langle \partial_xu_\infty(x_n),\dot{x}_n\right\rangle<\frac{1}{n}.
\end{equation}
It follows from Lemma \ref{apc} that there exists a positive constant $C$ independent of $x$ such that $|\partial_xu_\infty(x)|\leq C$ for $x\in \mathcal{D}$. Let $y_n:=\partial_xu_\infty(x_n)$.
    Extracting a subsequence if necessary, one obtain $(x_n,t_n,\dot{x}_n,y_n)\rightarrow(\bar{x},\bar{t},\dot{\bar{x}},\bar{y})$.  From (\ref{geq0}), we have
\[L(x_n,w(x_n,t_n),\dot{x}_n)-\left\langle y_n,\dot{x}_n\right\rangle\geq 0,\]which together with (\ref{555}) implies
\begin{equation}\label{777}
L(\bar{x},w(\bar{x},\bar{t}),\dot{\bar{x}})-\left\langle \bar{y},\dot{\bar{x}}\right\rangle=0.
\end{equation}Based on Lemma \ref{squr}, it follows from (\ref{666}) that for any $\xi\in T_{\bar{x}}M$, we have
\[L(\bar{x},w(\bar{x},\bar{t}),\xi)-\left\langle \bar{y},\xi\right\rangle\geq 0,\]which together with (\ref{777}) yields
\[\bar{y}=\frac{\partial L}{\partial \dot{x}}(\bar{x},w(\bar{x},\bar{t}),\dot{\bar{x}}).\]Hence,
\begin{equation}\label{888}
(\bar{x},w(\bar{x},\bar{t}),\dot{\bar{x}})\in \Delta.
 \end{equation}Since  $(x_n,w(x_n,t),\dot{x}_n)\in \Sigma(\delta)$, then
$(\bar{x},w(\bar{x},\bar{t}),\dot{\bar{x}})\in \Sigma(\frac{\delta}{2})$.
Moreover, there exists $\theta(\frac{\delta}{2})>0$ such that
\[\text{dist}((\bar{x},w(\bar{x},\bar{t}),\dot{\bar{x}}),\Delta)\geq \theta(\frac{\delta}{2}),\]
which is in contradiction with (\ref{888}). This verifies the claim.
\End

By Theorem \ref{one}, there exists a calibrated curve $\gm_w(t) : [0,s]\rightarrow M$ of $w(x,t)$ with $\gm_w(s)=x_0$.
In terms of (\ref{key}), it follows from (\ref{13}) that for any $\tau\in [0,s]$,
\begin{equation}\label{110}
H\left(\gm_w(\tau),w(\gm_w(\tau),\tau),\frac{\partial L}{\partial \dot{x}}(\gm_w(\tau),w(\gm_w(\tau),\tau),\dot{\gm}_w(\tau))\right)\geq\delta,
\end{equation}where $\frac{\partial L}{\partial \dot{x}}$ denotes the partial derivative of $L$ with respect to the third argument. By the Legendre transformation, we have $(\gm_w(\tau),w(\gm(\tau),\tau),\dot{\gm}_w(\tau))\in \Delta(\delta)$ for any $\tau\in [0,s]$.
It follows that for any $\tau\in [0,s]$, \[\text{dist}\left((\gm_w(\tau),w(\gm(\tau),\tau),\dot{\gm}_w(\tau)),\Delta\right)\geq\theta(\delta).\]
Let $\Theta$ be the set of $\gm(\tau)$ along which the directional derivative $\partial_{\gm_w(\tau)}u_\infty(\dot{\gm}_w(\tau))$ exists. For $\gm_w(\tau)\in \Theta$, we denote
\begin{equation}\label{120}
\widehat{L}(\gm_w(\tau)):= L(\gm_w(\tau),w(\gm_w(\tau),\tau),\dot{\gm}_w(\tau))-\partial_{\gm_w(\tau)}u_\infty(\dot{\gm}_w(\tau)).
\end{equation}
Since $u_\infty$ is a weak KAM solution, then it is locally semiconcave. By \cite{CS}, one can find a sequence $x^\tau_n\in \mathcal{D}$ with $x^\tau_n\rightarrow \gm_w(\tau)$ as $n\rightarrow\infty$ for a given $\tau\in [0,s]$ such that
\begin{equation}\label{77}
\partial_{\gm_w(\tau)}u_\infty(\dot{\gm}_w(\tau))\leq \langle \partial_{x}u_\infty(x^\tau_n),\dot{\gm}_w(\tau)\rangle+\frac{1}{n}.
\end{equation}
For $n$ large enough, it follows from (\ref{110}) that
\begin{equation}
H\left(x^\tau_n,w(x^\tau_n,\tau),\frac{\partial L}{\partial \dot{x}}(x^\tau_n,w(x^\tau_n,\tau),\dot{\gm}_w(\tau))\right)\geq\frac{\delta}{2},
\end{equation}which implies
\[\text{dist}\left((x^\tau_n,w(x^\tau_n,\tau),\dot{\gm}_w(\tau)),\Delta\right)\geq\theta(\frac{\delta}{2}).\]
Since $x^\tau_n\in \mathcal{D}$ and $(x^\tau_n,w(x^\tau_n,\tau),\dot{\gm}_w(\tau))\in \Sigma(\frac{\delta}{2})$, then there exists $\delta''>0$ independent of $\tau$  and $n$ such that for any $\tau\in [0,s]$,
\begin{equation}
\check{L}(x^\tau_n,w(x^\tau_n,\tau),\dot{\gm}_w(\tau))=L(x^\tau_n,w(x^\tau_n,\tau),\dot{\gm}_w(\tau))-\langle \partial_{x}u_\infty(x^\tau_n),\dot{\gm}_w(\tau)\rangle\geq \delta'',
\end{equation}which together with (\ref{120}) and (\ref{77}) implies
\begin{equation}
\widehat{L}(\gm_w(\tau))\geq \frac{\delta''}{2}.
\end{equation}
Moreover, we have
\begin{align}\label{18}
\int_0^s\widehat{L}(\gm_w(\tau))d\tau\geq \frac{\delta''}{2}s.
\end{align}
On the other hand, since $\gm_w$ is a calibrated curve of $w(x,t)$, then we have
\begin{align*}
\int_0^s\widehat{L}&(\gm_w(\tau))d\tau\\
&=\int_0^sL(\gm_w(\tau),w(\gm_w(\tau),\tau),\dot{\gm}_w(\tau))-\partial_{\gm_w(\tau)}u_\infty(\dot{\gm}_w(\tau))d\tau,\\
&\leq w(x_0,s)-w(\gm_w(0),0)-\left(u_\infty(\gm_w(s))-u_\infty(\gm_w(0))\right).
\end{align*}
By Proposition \ref{uniform bounds for semigroup}, there exists a positive constant $C$  independent of $s$ such that \[\int_0^s\widehat{L}(\gm_w(\tau))d\tau\leq C,\]
which is in contradiction with (\ref{18}) for $s$ large enough. Therefore, it yields that for $x\in \mathcal{D}$,
\begin{equation}\label{88}
H(x,\tilde{u}(x),\partial_x\tilde{u}(x))\leq 0.
\end{equation}
By virtue of a similar argument, one can obtain
\begin{equation}
H(x,\tilde{u}(x),\partial_x\tilde{u}(x))\geq 0,
\end{equation}which together with (\ref{88}) implies that for $x\in \mathcal{D}$, $H(x,\tilde{u}(x),\partial_x\tilde{u}(x))=0$.
This completes the proof of Lemma \ref{uu}.
\End

\subsection{Step 4}
Based on the preparations above, we will complete the proof of Theorem \ref{five}. Theorem \ref{six} can be concluded by Theorem \ref{five} directly.

\textbf{Proof of Theorem \ref{five}}:
Let $v(x,t):=T_t\tilde{u}(x)$, then $v(x,t)$ satisfies the following equation:
\begin{equation}\label{00}
\begin{cases}
\partial_tv(x,t)+H(x,v(x,t),\partial_xv(x,t))=0,\\
v(x,0)=\tilde{u}(x).
\end{cases}
\end{equation} Hence, $v(x,t)$ is locally Lipschitz on $M\times [0,\infty)$. Let $\mathcal{E}$ be the set of all differentiable points of $v$ on $M\times [0,\infty)$. Then $\mathcal{E}$ has full Lebesgue measure. Since $v(x,t)=T_{t+s}u^-(x)$, it follows from a similar argument as Lemma \ref{uu} that
 for $(x,t)\in \mathcal{E}$,
 \begin{equation}\label{xx}
H(x,v(x,t),\partial_x v(x,t)) = 0.
\end{equation}
%

For $(x,t)\in \mathcal{E}$, it follows from (\ref{00}) that
 \begin{equation}
\partial_t v(x,t)=-H(x,v(x,t),\partial_x v(x,t)).
\end{equation}By virtue of Fubini's Theorem, it follows from (\ref{xx}) that for almost every $x\in M$, $\partial_\tau v(x,\tau)= 0$ holds almost everywhere on $[0,t]$, where $\tau\in [0,t]$. Hence, we have
\begin{align*}
v(x,t)-v(x,0)=\int_0^t\partial_\tau v(x,\tau)d\tau= 0.
\end{align*}Since $v(x,t)$ is continuous, then for a given $(x,t)\in M\times [0,\infty)$, we have $v(x,t)=v(x,0)$. Combining with $v(x,t)=T_t\tilde{u}(x)$, it follows that
 $T_t\tilde{u}(x)= \tilde{u}(x)$ for any $t\geq 0$.


It follows that
 $\tilde{u}$ is a fixed point of $T_t$ for $t\geq 0$. By virtue of Lemma \ref{uisw}, we have $\tilde{u}$ is a weak KAM solution of (\ref{stahh}).
Moreover, using the non-expansiveness of $T_t$ again,  it follows that for $t>t_n$, we have
\begin{equation*}
\|T_t \phi - \tilde{u}\|_\infty = \|T_{t-t_n} \circ T_{t_n} \phi - T_{t-t_n}\tilde{u}\|_\infty \leq \|T_{t_n} \phi - \tilde{u}\|_\infty.
\end{equation*} Since $T_{t_n}\phi \rightarrow \tilde{u}$ as $t_n\rightarrow\infty$, we obtain
\[\lim_{t\rightarrow\infty}T_t\phi=\tilde{u},\]
where $\tilde{u}$ is a weak KAM solution of (\ref{stahh}).
This finishes the proof of Theorem \ref{five}.
\End

\begin{Remark}
Based on the uniqueness of the fixed point of $T_t$, we know that $u_\infty$ given by Lemma \ref{common fixed point existence} is the same as $\tilde{u}$ above.
\end{Remark}

 \vspace{2ex}
\noindent\textbf{Acknowledgement}
This work is partially under the support of National Natural Science Foundation of China (Grant No. 11171071,  11325103) and
National Basic Research Program of China (Grant No. 11171146).
X. Su is supported by both National
Natural Science Foundation of China (Grant No. 11301513) and ``the Fundamental Research Funds for the Central Universities".

{\sc Xifeng Su}

{\sc School of Mathematics, Beijing Normal University,
Beijing 100875,
China.}

 {\it E-mail address:} \texttt{xfsu@bnu.edu.cn}

\vspace{1em}

{\sc Lin Wang}

{\sc School of Mathematical Sciences, Fudan University,
Shanghai 200433,
China.}

 {\it E-mail address:} \texttt{linwang.math@gmail.com}

\vspace{1em}

{\sc Jun Yan}

{\sc School of Mathematical Sciences, Fudan University,
Shanghai 200433,
China.}

 {\it E-mail address:} \texttt{yanjun@fudan.edu.cn}

\end{document}